\documentclass[12pt]{article}
\usepackage{latexsym, amssymb}

\DeclareMathAlphabet\gothic{U}{euf}{m}{n}

\makeatletter
\def\eqnarray{\stepcounter{equation}\let\@currentlabel=\theequation
\global\@eqnswtrue
\tabskip\@centering\let\\=\@eqncr
$$\halign to \displaywidth\bgroup\hfil\global\@eqcnt\z@
  $\displaystyle\tabskip\z@{##}$&\global\@eqcnt\@ne
  \hfil$\displaystyle{{}##{}}$\hfil
  &\global\@eqcnt\tw@ $\displaystyle{##}$\hfil
  \tabskip\@centering&\llap{##}\tabskip\z@\cr}

\def\endeqnarray{\@@eqncr\egroup
      \global\advance\c@equation\m@ne$$\global\@ignoretrue}

\def\@yeqncr{\@ifnextchar [{\@xeqncr}{\@xeqncr[5pt]}}
\makeatother

\begin{document}
\bibliographystyle{tom}

\newtheorem{lemma}{Lemma}[section]
\newtheorem{thm}[lemma]{Theorem}
\newtheorem{cor}[lemma]{Corollary}
\newtheorem{voorb}[lemma]{Example}
\newtheorem{rem}[lemma]{Remark}
\newtheorem{prop}[lemma]{Proposition}
\newtheorem{stat}[lemma]{{\hspace{-5pt}}}
\newtheorem{obs}[lemma]{Observation}
\newtheorem{defin}[lemma]{Definition}

\newenvironment{remarkn}{\begin{rem} \rm}{\end{rem}}
\newenvironment{exam}{\begin{voorb} \rm}{\end{voorb}}
\newenvironment{defn}{\begin{defin} \rm}{\end{defin}}
\newenvironment{obsn}{\begin{obs} \rm}{\end{obs}}

\newenvironment{emphit}{\begin{itemize} }{\end{itemize}}

\newcommand{\gota}{\gothic{a}}
\newcommand{\gotb}{\gothic{b}}
\newcommand{\gotc}{\gothic{c}}
\newcommand{\gote}{\gothic{e}}
\newcommand{\gotf}{\gothic{f}}
\newcommand{\gotg}{\gothic{g}}
\newcommand{\gothh}{\gothic{h}}
\newcommand{\gotk}{\gothic{k}}
\newcommand{\gotm}{\gothic{m}}
\newcommand{\gotn}{\gothic{n}}
\newcommand{\gotp}{\gothic{p}}
\newcommand{\gotq}{\gothic{q}}
\newcommand{\gotr}{\gothic{r}}
\newcommand{\gots}{\gothic{s}}
\newcommand{\gotu}{\gothic{u}}
\newcommand{\gotv}{\gothic{v}}
\newcommand{\gotw}{\gothic{w}}
\newcommand{\gotz}{\gothic{z}}
\newcommand{\gotA}{\gothic{A}}
\newcommand{\gotB}{\gothic{B}}
\newcommand{\gotG}{\gothic{G}}
\newcommand{\gotL}{\gothic{L}}
\newcommand{\gotS}{\gothic{S}}
\newcommand{\gotT}{\gothic{T}}

\newcommand{\mn}{\marginpar{\hspace{1cm}*} }
\newcommand{\mnn}{\marginpar{\hspace{1cm}**} }

\newcommand{\mnq}{\marginpar{\hspace{1cm}*???} }
\newcommand{\mnnq}{\marginpar{\hspace{1cm}**???} }

\newcounter{teller}
\renewcommand{\theteller}{\Roman{teller}}
\newenvironment{tabel}{\begin{list}%
{\rm \bf \Roman{teller}.\hfill}{\usecounter{teller} \leftmargin=1.1cm
\labelwidth=1.1cm \labelsep=0cm \parsep=0cm}
                      }{\end{list}}

\newcounter{tellerr}
\renewcommand{\thetellerr}{(\roman{tellerr})}
\newenvironment{subtabel}{\begin{list}%
{\rm  (\roman{tellerr})\hfill}{\usecounter{tellerr} \leftmargin=1.1cm
\labelwidth=1.1cm \labelsep=0cm \parsep=0cm}
                         }{\end{list}}
\newenvironment{ssubtabel}{\begin{list}%
{\rm  (\roman{tellerr})\hfill}{\usecounter{tellerr} \leftmargin=1.1cm
\labelwidth=1.1cm \labelsep=0cm \parsep=0cm \topsep=1.5mm}
                         }{\end{list}}

\newcommand{\Ni}{{\bf N}}
\newcommand{\Ri}{{\bf R}}
\newcommand{\Ci}{{\bf C}}
\newcommand{\Si}{{\bf S}}
\newcommand{\Ti}{{\bf T}}
\newcommand{\Zi}{{\bf Z}}
\newcommand{\Fi}{{\bf F}}

\newcommand{\Epsilon}{{\rm E}}

\newcommand{\proof}{\mbox{\bf Proof} \hspace{5pt}} 
\newcommand{\remark}{\mbox{\bf Remark} \hspace{5pt}}
\newcommand{\ruimte}{\vskip10.0pt plus 4.0pt minus 6.0pt}

\newcommand{\simh}{{\stackrel{{\rm cap}}{\sim}}}
\newcommand{\ad}{{\mathop{\rm ad}}}
\newcommand{\Ad}{{\mathop{\rm Ad}}}
\newcommand{\Aut}{\mathop{\rm Aut}}
\newcommand{\arccot}{\mathop{\rm arccot}}
\newcommand{\capp}{{\mathop{\rm cap}}}
\newcommand{\rcapp}{{\mathop{\rm rcap}}}
\newcommand{\Capp}{{\mathop{\rm Cap}}}
\newcommand{\diam}{\mathop{\rm diam}}
\newcommand{\divv}{\mathop{\rm div}}
\newcommand{\dist}{\mathop{\rm dist}}
\newcommand{\codim}{\mathop{\rm codim}}
\newcommand{\RRe}{\mathop{\rm Re}}
\newcommand{\IIm}{\mathop{\rm Im}}
\newcommand{\Tr}{{\mathop{\rm Tr}}}
\newcommand{\Vol}{{\mathop{\rm Vol}}}
\newcommand{\card}{{\mathop{\rm card}}}
\newcommand{\supp}{\mathop{\rm supp}}
\newcommand{\sgn}{\mathop{\rm sgn}}
\newcommand{\essinf}{\mathop{\rm ess\,inf}}
\newcommand{\esssup}{\mathop{\rm ess\,sup}}
\newcommand{\Int}{\mathop{\rm Int}}
\newcommand{\Leibniz}{\mathop{\rm Leibniz}}
\newcommand{\lcm}{\mathop{\rm lcm}}
\newcommand{\loc}{{\rm loc}}

\newcommand{\mod}{\mathop{\rm mod}}
\newcommand{\spann}{\mathop{\rm span}}
\newcommand{\one}{1\hspace{-4.5pt}1}

\newcommand{\DWR}{}

\hyphenation{groups}
\hyphenation{unitary}

\newcommand{\tfrac}[2]{{\textstyle \frac{#1}{#2}}}

\newcommand{\cb}{{\cal B}}
\newcommand{\cc}{{\cal C}}
\newcommand{\cd}{{\cal D}}
\newcommand{\ce}{{\cal E}}
\newcommand{\cf}{{\cal F}}
\newcommand{\ch}{{\cal H}}
\newcommand{\ci}{{\cal I}}
\newcommand{\ck}{{\cal K}}
\newcommand{\cl}{{\cal L}}
\newcommand{\cm}{{\cal M}}
\newcommand{\cn}{{\cal N}}
\newcommand{\co}{{\cal O}}
\newcommand{\cs}{{\cal S}}
\newcommand{\ct}{{\cal T}}
\newcommand{\cx}{{\cal X}}
\newcommand{\cy}{{\cal Y}}
\newcommand{\cz}{{\cal Z}}

\newcommand{\wtozp}{W^{1,2}\raisebox{10pt}[0pt][0pt]{\makebox[0pt]{\hspace{-34pt}$\scriptstyle\circ$}}}
\newlength{\hightcharacter}
\newlength{\widthcharacter}
\newcommand{\covsup}[1]{\settowidth{\widthcharacter}{$#1$}\addtolength{\widthcharacter}{-0.15em}\settoheight{\hightcharacter}{$#1$}\addtolength{\hightcharacter}{0.1ex}#1\raisebox{\hightcharacter}[0pt][0pt]{\makebox[0pt]{\hspace{-\widthcharacter}$\scriptstyle\circ$}}}
\newcommand{\cov}[1]{\settowidth{\widthcharacter}{$#1$}\addtolength{\widthcharacter}{-0.15em}\settoheight{\hightcharacter}{$#1$}\addtolength{\hightcharacter}{0.1ex}#1\raisebox{\hightcharacter}{\makebox[0pt]{\hspace{-\widthcharacter}$\scriptstyle\circ$}}}
\newcommand{\scov}[1]{\settowidth{\widthcharacter}{$#1$}\addtolength{\widthcharacter}{-0.15em}\settoheight{\hightcharacter}{$#1$}\addtolength{\hightcharacter}{0.1ex}#1\raisebox{0.7\hightcharacter}{\makebox[0pt]{\hspace{-\widthcharacter}$\scriptstyle\circ$}}}

\newpage

 \thispagestyle{empty}
  
 \begin{center}
 \vspace*{-1.0cm}
\vspace*{1.5cm}

{\Large{\bf On self-adjointness  of }}\\[3mm]
{\Large{\bf symmetric diffusion operators }}  \\[5mm]
\large Derek W. Robinson$^\dag$ \\[1mm]

\normalsize{8th November 2019}\\[1mm]
\end{center}

\vspace{+5mm}

\begin{list}{}{\leftmargin=1.7cm \rightmargin=1.7cm \listparindent=15mm 
   \parsep=0pt}
   \item
{\bf Abstract} $\;$ 
Let $\Omega$ be a domain in $\Ri^d$ with  boundary $\Gamma$ and let $d_\Gamma$ denote the Euclidean distance to  $\Gamma$.
Further let $H=-\divv(C\nabla)$ where  $C=(\,c_{kl}\,)>0$ with $c_{kl}=c_{lk}$ are real, bounded, Lipschitz continuous functions and 
$D(H)=C_c^\infty(\Omega)$.
Assume also that there is a $\delta\geq0$ such that  $\|C/d_\Gamma^{\,\delta}-aI\|\to 0$ as $d_\Gamma\to0$ with $\delta\geq0$
where $a$ is a bounded Lipschitz continuous function with $a\geq\mu>0$ on a boundary layer $\Gamma_{\!\!r}=\{x\in\Omega: d_\Gamma(x)<r\}$.
Finally we require  $|(\divv C).(\nabla d_\Gamma)|d_\Gamma^{\,-\delta+1}$ to be bounded  on~$\Gamma_{\!\!r}$.
Then we prove that if $\Omega$ is  a $C^2$-domain, or if $\Omega=\Ri^d\backslash S$ where $S$ is a countable set of positively separated points, or if 
$\Omega=\Ri^d\backslash \overline \Pi$ with $\Pi$ a convex set whose  boundary has 
Hausdorff dimension $d_H\in \{1,\ldots, d-1\}$ then  the condition $\delta>2-(d-d_H)/2$ is sufficient for $H$ to be 
essentially self-adjoint as an operator on $L_2(\Omega)$.
In particular $\delta>3/2$ suffices for $C^2$-domains.
Finally we prove that  $\delta\geq 3/2$ is  necessary in the $C^2$-case.
\end{list}

\vfill



\vspace{0.5cm}

\noindent
\begin{tabular}{@{}cl@{\hspace{10mm}}cl}
$ {}^\dag\hspace{-5mm}$&   Mathematical Sciences Institute (CMA)    &  {} &{}\\
  &Australian National University& & {}\\
&Canberra, ACT 0200 && {} \\
  & Australia && {} \\
  &derek.robinson@anu.edu.au
 & &{}\\
\end{tabular}

\newpage

\setcounter{page}{1}

\section{Introduction}\label{S1}
Our intention is to analyze the $L_2$-uniqueness of symmetric diffusion processes on a domain $\Omega$ of the Euclidean space $\Ri^d$
with boundary $\Gamma$.
The problem can either be expressed as uniqueness of weak solutions of a diffusion equation
\[
\partial\varphi_t/\partial t+H\varphi_t=0
\]
on $L_2(\Omega)$ or, equivalently, as essential self-adjointness of the diffusion operator  $H=-\divv(C\,\nabla)$ on $C_c^\infty(\Omega)$.
In the standard theory of strongly elliptic operators uniqueness is ensured by the specification of  conditions at the boundary.
If, however, the coefficient matrix $C$ is degenerate at the boundary then uniqueness is determined by the properties of the diffusion in a neighbourhood of $\Gamma$.
In both cases, however, the existence of a solution to the diffusion equation follows by specifying Dirichlet boundary conditions or,
in operator terms, by constructing the self-adjoint Friedrichs' extension $H_{\!F}$ of $H$.
So in  the degenerate case  the uniqueness problem consists of relating the boundary behaviour to the uniqueness of the Dirichlet solution or the Friedrichs' extension.

Some guidance to  $L_2$-uniqueness is given by the better understood problem of  $L_1$-uniqueness.
The Friedrichs' extension $H_{\!F}$ generates a submarkovian semigroup $S^F$ on $L_2(\Omega)$, i.e.\ a semigroup which is both positive
and $L_\infty$-contractive, and this leads to a natural probabilistic interpretation.
Then $L_1$-uniqueness is   equivalent to $S^F$ conserving probability by Theorem~2.2 in \cite{Dav14} and this  in turn is equivalent to $H_{\!F}$ being the 
unique self-adjoint extension of $H$ which generates a submarkovian semigroup (see \cite{Ebe}, Corollary~3.4,  and  \cite{RSi4}, Theorem~1.3).
Thus the $L_1$-problem reduces to an $L_2$-problem, Markov uniqueness,   which turns out to be considerably simpler than essential self-adjointness.
The $L_1$-uniqueness is a problem of quadratic forms  but the $L_2$-uniqueness is a genuine operator problem.
Background to these problems in various contexts, e.g.\ diffusion on Riemannian manifolds, stochastic differential equations,
Dirichlet forms, etc., for various kinds of diffusion, e.g.\ with or without drift, stochastically complete or incomplete, can be found  in \cite{Dav14} \cite{FOT} \cite{Ebe}
\cite{RSi4} or the more recent \cite{NeN}.
The introduction to the latter paper gives a broad description of the problems of interest.

\smallskip

We assume throughout that $C=(\,c_{kl}\,)>0$ with $c_{kl}=c_{lk}$ real, bounded, Lipschitz continuous functions over $\Omega$
which are degenerate at the boundary.
The rate of degeneracy is measured in terms of $d_\Gamma$, the Euclidean distance to the boundary.
In particular we allow the coefficients $c_{kl}$ and their derivatives  to behave like $d_\Gamma^{\,\delta}$ and $d_\Gamma^{\,\delta-1}$,
respectively,  with $\delta\geq 0$ as $d_\Gamma\to0$.
First we establish under a very mild degeneracy assumption that the $L_2$-uniqueness problem is reduced to understanding the properties of $H_{\!F}$ in a thin boundary layer $\Gamma_{\!\!r}=\{x\in\Omega: d_\Gamma(x)<r\}$.
The term boundary layer is adopted from the theory of laminar flow where the rate of 
diffusion adjacent to the boundary is governed by the distance from the boundary.
It appears appropriate in the current context although the sets $\Gamma_{\!\!r}$ are variously described as 
tubular neighbourhoods,  inner $r$-neighbourhoods, parallel bodies, etc.\ in other  mathematical areas.
Secondly, as a consequence of this separation, $L_2$-uniqueness is reduced to identifying the Friedrichs' extension with the closure $\overline H$ of $H$ as operators on $L_2(\Gamma_{\!\!r})$.
Since the second-order operator $H$ is `comparable' to the multiplication operator $d_\Gamma^{\,\delta-2}$ near the boundary one would expect
that $D_r=\{\varphi\in D(H_{\!F}): \supp\varphi \subseteq \overline\Gamma_{\!\!r}\}\subset D(d_\Gamma^{\,\delta-2})$ for all small~$r$.
This is a  form of the well known Rellich inequality and presents a principal  technical problem in deriving sufficiency conditions for uniqueness.
It is clearly valid if $\delta\geq2$  and indeed $L_2$-uniqueness then follows from a standard criterion for self-adjointness of elliptic operators
developed in the 1960s and 70s which can be found in \cite{Dav14}, Section~3, together with background references (see also  \cite{ERS5}, Corollary~4.12).
It is notable that the $\delta\geq 2$ result is independent of the geometry of the boundary and requires no particular constraints on the derivatives of the coefficients $c_{kl}$.
Understanding the situation for smaller values of $\delta$ is much more sensitive to  properties of the boundary and the coefficients.

The key criterion for $L_1$-uniqueness  which includes values of $\delta<2$  is  given by 
\begin{equation}
\delta\geq 2-(d-d_H)
\label{eun1.2}
\end{equation}
where $d_H$ is the Hausdorff dimension of the boundary $\Gamma$.
The  sufficiency of the condition was first proved  for a large variety of domains in \cite{RSi4} 
(see also \cite{ERS5} Theorems~4.10 and 4.14).
Subsequently,  the proof was extended to all Ahlfors regular domains by Proposition~2.4  in \cite{LR}.
This includes domains with rough, smooth, fractal or uniformly disconnected boundaries (see  \cite{Hei} or   \cite{Semm}  for details on the Ahlfors property). 
 In addition Proposition~3.6 of \cite{LR} established that  (\ref{eun1.2}) is a necessary condition for  $L_1$-uniqueness  for almost all Ahlfors regular domains.
The latter result  is not quite universal  because of various possible boundary pathologies such as cuspoidal points or antennae.
Condition~(\ref{eun1.2}) is directly related to the local integrability of $d_\Gamma^{\,\delta-2}$ on ${\overline\Gamma}_{\!\!r}$.
This depends not only on the order of growth,  $2-\delta$, but also the measure of subsets of the boundary layer.
  The analysis of \cite{RSi4} \cite{ERS5} and \cite{LR} indicates that the condition analogous to (\ref{eun1.2}) 
  for $L_2$-uniqueness should be  
\begin{equation}
\delta\geq 2-(d-d_H)/2
\label{eun1.3}
\end{equation}
and $\delta\geq 2-(d-d_H)/p$ for $L_p$-uniqueness.
The latter condition is related to the local $L_p$-integrability of  $d_\Gamma^{\,\delta-2}$ on the boundary layers.
But this is almost certainly an oversimplification.
It is very likely that more detailed geometric characteristics  of the boundary than its Hausdorff dimension are important.
The $L_1$-analysis also indicates that necessary conditions for uniqueness could be more sensitive to boundary properties than sufficient conditions.
Most of the results in the sequel require some smoothness of the boundary.

First we establish that if   $\Omega$ is a $C^2$-domain 
then $H$ is essentially self-adjoint if $\delta>3/2$.
Since the $C^2$-property ensures that $d_H=d-1$ this is 
 in agreement  with (\ref{eun1.3}).
Note that the $C^2$-property is usually defined locally (see,  for example,  Section~6.2 of \cite{GT}) and the definition guarantees for bounded domains  that the principal curvatures of the boundary are  uniformly bounded.
Nevertheless the definition  extends in a natural manner to unbounded domains  retaining the boundedness (see \cite{For}, Appendix~B).
Our results encompass the unbounded case.
Moreover, if $\Omega$ is a $C^2$-domain  the exterior $\Ri^d\backslash\overline\Omega$ is also a $C^2$-domain and since   both domains have a common boundary our conclusion is  equally valid.
Secondly, we prove that if $\Omega=\Ri^d\backslash\{0\}$ and $\delta>2-d/2$ then $H$ is essentially self-adjoint.
 This  second case can be considered a degenerate example of the $C^2$-exterior domain.
Thirdly, we derive analogous conclusions for a variety of intermediate situations.
 If $\Pi$ is a convex  $C^2$-domain in $\Ri^s$, with $s\in \{1,\ldots, d-1\}$, and $\Omega=\Ri^d\backslash\overline\Pi$ then the boundary $\Gamma$ of $\Omega$ is equal to $\overline\Pi$ and $d_H=s$.
 In this case we establish that $\delta>2-(d-d_H)/2$ is sufficient for self-adjointness.
Since each of these results  is deduced from estimates on an arbitrarily thin boundary layer and since the boundary layers decompose
into disjoint components if $\Gamma$ consists of positively separated components, e.g.\ if $\Omega$ is an annulus, the special cases
can be combined to give many more general conclusions.
For example, if $\Omega=\Ri^d\backslash\Zi^d$ then $\delta>2-d/2$ is again sufficient for self-adjointness.
A precise formulation of these sufficiency results is given in Theorem~\ref{tun5.1} and the discussion following the theorem.
Our sufficiency results partially overlap with the conclusions of \cite{NeN} although our methods and approaches are quite different.

Finally  we establish that the condition $\delta\geq 3/2$ is necessary for self-adjointness  for all  $C^2$-domains.
The precise result is given in Theorem~\ref{tun5.11} which does require slightly more stringent bounds on the derivatives of the coefficients of $H$ than those used  in Theorem~\ref{tun5.1}.
Nevertheless the combined results  establish that (\ref{eun1.3}) is both necessary and sufficient for $L_2$-uniqueness on $C^2$-domains with the exception of the sufficiency of 
the critical value $\delta=3/2$.
The proof of necessity in the $C^2$-case is facilitated by the fact that the boundary has codimension $1$.
Even the seemingly simple case of $\Ri^d\backslash\{0\}$ is complicated for general operators because the codimension is large.

We begin in  Section~\ref{S2} with a more precise definition of the operators we subsequently  examine and then we discuss successively 
the boundary localization technique,  the construction of approximate identities, the differentiability properties of $d_\Gamma$ and finally Hardy--Rellich inequalities near the boundary.
In Section~\ref{S3} we apply these techniques to establish the sufficiency results for $L_2$-uniqueness.
Then we recall some elementary results on Markov uniqueness and show how these can be utilized to obtain the necessary condition for $C^2$-domains.
We conclude with  comments on open problems in Section~\ref{S4}.

\section{Localization and approximation}\label{S2}

The elliptic operator $H=-\divv(C\nabla)$ is  initially  defined on the domain
$D(H)=C_c^\infty(\Omega)$.
The matrix $C=(\,c_{kl}\,)$ of coefficients is real, symmetric and strictly positive on $\Omega$.
Further the $c_{kl}$ are bounded, Lipschitz continuous functions.
It follows that $H$ is symmetric and consequently closable with respect to the graph norm $\|\varphi\|_{D(H)}=(\|H\varphi\|_2^2+\|\varphi\|_2^2)^{1/2}$.
Since  the coefficients $c_{kl}$ are bounded $\nu=\sup_{x\in\Omega}\|C(x)\|<\infty$, where $\|\cdot\|$ denotes the matrix norm,
and  $C\leq \nu \,I$.
Moreover, it follows from the strict positivity that for each compact subset $K$ of $\Omega$ there is a $\mu_{\!K}>0$ such that $C\geq \mu_{\!K} I$.
Thus $C$ is locally strongly elliptic. 
Hence, by elliptic regularity, the domain of the adjoint $H^*$ of  $H$ is contained in  $W^{2,2}_{\rm loc}(\Omega)$.
These properties will be assumed throughout.

The behaviour of the coefficients at the boundary $\Gamma$ is specified as follows.
We assume there are an $r_0\in\langle0,1]$, an $a\in W^{1,\infty}({\overline\Gamma_{\!\!r_0}})$  and $\lambda, \mu>0$ such that 
$\lambda\geq a\geq \mu$ and 
 \begin{equation}
\textstyle{\inf_{r\in\langle0,r_0]}}\;\textstyle{\sup_{x\in\Gamma_{\!\!r}}}\|C(x)d_\Gamma(x)^{-\delta}-a(x)I\|=0
\label{euns3.2}
 \end{equation}
 where $\delta\geq0$.
In addition
\begin{equation}
\textstyle{\sup_{x\in\Gamma_{\!\!r_0}}}\Big(|(\divv C).(\nabla\!d_\Gamma)(x)|\,d_\Gamma(x)^{-\delta+1}\Big)<\infty
\label{eun3.21}
 \end{equation}
where $\divv C$ is the vector with components  $(\divv C)_l=\sum^d_{k=1}\partial_kc_{kl}$. 
These conditions are a statement that $C$ is `comparable' to $ad_\Gamma^{\,\delta}$ at the boundary.
The  parameter $\delta$ clearly determines the order of degeneracy at the boundary and the function $a$ is a measure of the boundary  profile of the coefficients. 
Condition (\ref{eun3.21}) can be interpreted as a bound on the derivatives of the coefficients  in the direction normal to the boundary if this makes sense.
These conditions will be used in the subsequent derivation of sufficiency conditions for self-adjointness but we will require a slight strengthening of (\ref{eun3.21}) in the  
discussion of necessary conditions.
The simplest illustration of these conditions is given by multiples of the identity, e.g.\  $C=ad_\Gamma^{\,\delta}I$ on $\Gamma_{\!\!r}$.
Then $Cd_\Gamma^{\,-\delta}=a I$, $(\divv (Cd_\Gamma^{\,-\delta})).(\nabla\!d_\Gamma)=(\nabla\! a).(\nabla\!d_\Gamma)$ and 
$|(\divv C).(\nabla\!d_\Gamma)|\,d_\Gamma^{\,-\delta+1}\leq \delta a +|(\nabla \!a).(\nabla\!d_\Gamma)|\, d_\Gamma$.

\smallskip
The notion of comparability can be made more precise by noting
that  for each $r\in\langle0,r_0]$  there are $\sigma_{\!r}, \tau_{\!r}>0$ such that 
\begin{equation}
\sigma_{\!r}(a d_\Gamma^{\,\delta})(x) I\leq C(x)\leq \tau_{\!r} (a d_\Gamma^{\,\delta})(x) I
\label{euns2.2}
\end{equation}
for all $x\in \Gamma_{\!\!r}$.
Moreover, it follows from (\ref{eun3.21}) there is also a $\rho_r>0$ such that 
\begin{equation}
|((\divv C).(\nabla\!d_\Gamma))(x)|\leq \rho_r(a d_\Gamma^{\,\delta-1})(x)
\label{euns2.2.1}
\end{equation}
for all $x\in \Gamma_{\!\!r}$.
The earlier discussions of Markov uniqueness  in \cite{RSi4} and \cite{LR} were based on the conditions (\ref{euns2.2}) and the derivative bounds (\ref{euns2.2.1}) were unnecessary.
They will, however, play a vital role in the sequel.

\smallskip

Next let $h$ denote the  positive bilinear form
 associated with $H$ on $L_2(\Omega)$, i.e.\ 
\[
h(\psi, \varphi)
=(\psi, H\varphi)= \sum^d_{k,l=1}(\partial_k\psi, c_{kl} \, \partial_l\varphi)
\]
for all $\psi, \varphi\in C_c^\infty(\Omega)$ and set $h(\varphi)=h(\varphi, \varphi)$.
Since $H$ is positive-definite  and symmetric  the form $h$ is   
closable with respect to the graph norm $\|\varphi\|_{D(h)}=(h(\varphi)+\|\varphi\|_2^2)^{1/2}$.
For economy of notation we now use $H$ and $h$ to denote the closures of the  operator and  form, respectively. 
The closed form  $h$  is  a Dirichlet form \cite{BH} \cite{FOT}.
Therefore there   is a positive self-adjoint extension $H_{\!F}$, the  Friedrichs' extension of $H$, associated with $h$.
In fact $D(h)=D(H_{\!F}^{1/2})$ and $H_{\!F}$ generates a submarkovian semigroup.
Since  $H\subseteq H_{\!F}$  it follows that  $H$ is self-adjoint  if and only if $H_{\!F}\subseteq H$. 
This is the criterion we  use in the next section to establish self-adjointness.

\smallskip

The Dirichlet form $h$ has a {\it carr\'e du champ}, a positive continuous bilinear form $\psi,\varphi\in D(h)\mapsto \Gamma(\psi,\varphi)\in L_1(\Omega)$, such that 
$h(\psi,\varphi)=\|\Gamma(\psi,\varphi)\|_1$ (see, for example, \cite{BH} Section~I.4).
Moreover, $\Gamma$ satisfies the Cauchy--Schwarz inequality, i.e.\ 
$\|\Gamma(\psi, \varphi)\|_1^2\leq \|\Gamma(\psi)\|_1\, \|\Gamma(\varphi)\|_1$ with $\Gamma(\psi)=\Gamma(\psi,\psi)$.
If $\psi, \varphi\in C_c^\infty(\Omega)$ then
\[
\Gamma(\psi,\varphi)=\sum^d_{k,l=1}c_{kl} (\partial_k\psi) (\partial_l\varphi)
\]
and one can use this identification to define $\Gamma$ as a function over $W^{1,\infty}(\Omega)\times C_c^\infty(\Omega)$.
But then
\[
|(\psi, \Gamma(\chi,\varphi))|\leq h(\varphi)^{1/2}\,\|\Gamma(\chi)\|_\infty^{1/2}\,\|\psi\|_2
\]
for all $\psi,\varphi\in C_c^\infty(\Omega)$ and $\chi\in W^{1,\infty}(\Omega)$.
Therefore $\Gamma$ extends to an $L_2$-function on $W^{1,\infty}(\Omega)\times D(h)$ with 
$\|\Gamma(\chi,\varphi)\|_2^2\leq \|\Gamma(\chi)\|_\infty\,\|\Gamma(\varphi)\|_1$.

\subsection{Boundary localization}\label{S2-2}

In this subsection we establish that the criterion $H=H_{\!F}$ for self-adjointness of the operator $H$
can be reduced to a problem of boundary behaviour by localization arguments.
In particular if $\Omega_t=\{x\in\Omega: d_\Gamma(x)>t\}$ then  one can deduce by localization that $H_{\!F}\varphi=H\varphi$ for all $\varphi\in D(H_{\!F})$
 with  $\supp\varphi\subset \Omega_t$ for some $t>0$.
 Conditions~(\ref{euns2.2}) and (\ref{euns2.2.1}) are unnecessary for this conclusion. They can be replaced by the assumption
 that there is a $\mu_u>0$ such that $C(x)\geq \mu_uI$ for all $x\in\Omega_u\backslash\Omega_t$ for  small $t>u>0$.
 
 The proof begins with the following proposition.

\begin{prop}\label{pun2.1}
\begin{tabel}
\item\label{pun2.1-1} 
\hspace{-7mm} If $\chi\in W^{1,\infty}(\Omega)$ and $\varphi\in D(h)$ then $\chi\varphi\in D(h)$ and 
\begin{equation}
h(\chi\varphi)\leq 2\,\|\chi\|_\infty^2\, h(\varphi)+2\,\|\Gamma(\chi)\|_\infty\|\varphi\|_2^2\;.
\label{eun2.30}
\end{equation}
\item\label{pun2.1-2} 
\hspace{-5mm} If  $\chi\in W^{2,\infty}(\Omega)$ and $\varphi\in D(H_{\!F})$ then $\chi\varphi\in D(H_{\!F})$ and
\begin{equation}
H_{\!F}(\chi\varphi)=\chi (H_{\!F}\varphi)+(H\chi) \varphi-2\,\Gamma(\chi,\varphi)
\label{eun2.3}
\end{equation}
where $H\chi =-\divv(C\nabla\chi)$.
\end{tabel}
\end{prop}
\proof\ The proposition is a generalization  of a result for  $\Omega=\Ri^d$ given by  Lemma~4.3 of \cite{ERS5}.
The first statement is a corollary of a general  property of  local Dirichlet forms.
If $\chi\in W^{1,\infty}(\Omega)$ it follows that $\chi\in D(h)_{\rm loc}\cap L_\infty(\Omega)$ and $\|\Gamma(\chi)\|_\infty<\infty$. 
Therefore if $\varphi\in D(h)$ then  $\chi\varphi\in D(h)$ and  (\ref{eun2.30}) is valid by  Lemma~3.4 of \cite{ERSZ2}.

\smallskip
Next if $\chi\in W^{2,\infty}(\Omega)$
and $\psi, \varphi\in C_c^\infty(\Omega)$ then
\begin{equation}
h(\psi, \chi\varphi)=h(\chi\psi, \varphi)+(\psi,(H\chi)\,\varphi)
-2\,(\psi,\Gamma(\chi,\varphi))
\label{eun2.4}
\end{equation}
where $H\chi=-\divv(C\nabla\chi)\in L_\infty(\Omega)$.
Hence 
\begin{equation}
|h(\psi, \chi\varphi)|\leq|h(\chi\psi, \varphi)|+\Big(\|H\chi\|_\infty \,\|\varphi\|_2
+2\, \|\Gamma(\chi)\|_\infty^{1/2} \,h(\varphi)^{1/2}\Big)\,\|\psi\|_2
\label{eun2.5}
\end{equation}
and, by continuity,  (\ref{eun2.5}) extends to all $\psi,\varphi\in D(h)$.
But if $\varphi\in D(H_{\!F})$ then 
\[
|h(\chi\psi, \varphi)|=|(\chi\psi, H_{\!F}\varphi)|\leq \|H_{\!F}\varphi\|_2\,\|\chi\|_\infty\,\|\psi\|_2
\]
for all $\psi\in D(h)$.
Thus  it follows from (\ref{eun2.5}) that  there is an $a>0$ such that $|h(\psi, \chi\varphi)|\leq a\,\|\psi\|_2$
for all $\psi\in D(h)$.
Therefore $\chi\varphi\in D(H_{\!F})$.
Finally (\ref{eun2.3}) follows from  the extension of (\ref{eun2.4}) to $\psi\in D(h)$ and $\varphi\in D(H_{\!F})\subseteq D(h)$.
\hfill$\Box$

\bigskip

In order to exploit the strict positivity of the coefficient matrix $C$
we introduce a family of      extensions $ H_{\!u}$ of $H$ which act on $L_2(\Ri^d)$.
First fix $u,t>0$ with $u<t<1$.
Then $\Omega_t\subset \Omega_u$ and  it follows from the boundary degeneracy assumption (\ref{euns2.2}) that there is a $\mu_u>0$ such that $C\geq \mu_u I$ on  $\Omega_u\backslash\Omega_t$, e.g.\ one can take $\mu_u=\mu\,\sigma_{\!u}u^\delta$.
Secondly, choose $\xi\in W^{2,\infty}(\Ri^d)$ such that $0\leq \xi\leq 1$, $\supp\xi\subset \Omega_u$ and $\xi=1$ on $\Omega_t$.
Thirdly, introduce the coefficient matrix $ C_u=\xi\,C+(1-\xi)\,\mu_u I$ on $\Ri^d$.
Finally let $H_{\!u}$ denote  the symmetric operator on $C^\infty_c(\Ri^d)$ corresponding to the coefficient matrix $C_u$.
Specifically, $H_{\!u}=-\divv(C_u\nabla)$.
The operator is again closable and, for simplicity, we use $H_{\!u}$ to denote the  closure  and $h_u$ the corresponding
Dirichlet form.
Thus  $h_u$   is  the closure of the form 
$\varphi\in C_c^\infty(\Ri^d)\mapsto  h_u(\varphi)=(\varphi,  H_{\!u}\varphi)$.
It follows by construction that  $ H_{\!u}\varphi=H\varphi$ and $h_u(\varphi)=h(\varphi)$ for all $\varphi\in C_c^\infty(\Omega_t)$.

The latter  properties extend to the closed forms and operators by the following proposition.

\begin{prop}\label{pun2.11}
Let  $\supp\varphi\subset \Omega_r$  with $r>t>u$.  

\begin{tabel}
\item\label{pun2.11-1}
$\;\varphi\in D(h)$ if and  only if $\varphi\in D(h_u)$  and then $h(\varphi)=h_u(\varphi)$.

\item\label{pun2.11-2}
$\;\varphi\in D(H)$ if and only $\varphi\in D(H_{\!u})$ and then $H\varphi=H_u\varphi$.

\item\label{pun2.11-3}
$\:\varphi\in D(H_{\!F})$ if and only if $\varphi\in D((H_{\!u})_{\!F})$ and then $H_{\!F}\varphi=(H_{\!u})_{\!F}\varphi$.
\end{tabel}
\end{prop}
\proof\
\ref{pun2.11-1}.$\;$
Assume $\varphi\in D(h)$.
One can choose a sequence $\varphi_n\in C_c^\infty(\Omega)$ such that $\|\varphi_n-\varphi\|_{D(h)}\to0$ as $n\to\infty$.
Since $r>t$ one can also choose $\chi\in W^{2,\infty}(\Omega_t)$ with $0\leq \chi\leq 1$ and $\chi=1$ on $\Omega_r$.
Then $\chi\varphi=\varphi$ and $\|\chi\varphi_n-\varphi\|_{D(h)}=\|\chi(\varphi_n-\varphi)\|_{D(h)}\to0$ as $n\to\infty$ as a simple consequence of (\ref{eun2.30}).
Since $h_u(\chi(\varphi_n-\varphi_m))=h(\chi(\varphi_n-\varphi_m))$, by the above discussion, it follows that $\varphi\in D( h_u)$,
$\|\chi\varphi_n-\varphi\|_{D(h)}\to0$ and $h_u(\varphi)=h(\varphi)$.
The converse is proved similarly.

\smallskip

\ref{pun2.11-2}.$\;$
 The proof is similar with graph norms $\|\cdot\|_{D(h)}$  of the form replaced by the graph norms $\|\cdot\|_{D(H)}$ of the operators.

\smallskip

\ref{pun2.11-3}.$\;$
If $\varphi\in D(H_{\!F})\subset D(h) $ with $\supp\varphi\subset \Omega_r$ 
then $\varphi\in D(h_u)$, by the foregoing.
Now  if $\psi\in D(h)$ then $\chi\psi\in D(h)$ and $\supp\chi\psi\subset \Omega_t$.
Therefore $\chi\psi\in D(h_u)$ and $h(\chi\psi,\varphi)=h_u(\chi\psi,\varphi)$.
But $h$ and $ h_u$ are both local Dirichlet forms in the sense of \cite{BH} (or strongly local in the terminology of \cite{FOT}).
Therefore
\[
 h_u(\psi,\varphi)= h_u(\chi\psi, \varphi)+h_u((\one-\chi)\psi, \varphi)=h_u(\chi\psi, \varphi)=h(\chi\psi,\varphi)
\]
since 
 $h_u((\one-\chi)\psi, \varphi)=0$   by locality of $h_u$.
Hence $ h_u(\psi,\varphi)= h(\chi\psi, \varphi)$.
Then, however, 
\[
 |h_u(\psi,\varphi)|= |h(\chi\psi, \varphi)|= |(\chi\psi, H_{\!F}\varphi)|\leq \|\psi\|_2\,\|H_{\!F}\varphi\|_2
 \;.
 \]
Hence $\varphi\in D((H_{\!u})_{\!F})$.
Interchanging $H$ and $H_{\!u}$ in this argument one concludes that if  $\varphi\in D((H_{\!u})_{\!F})$ with 
$\supp\varphi\subset \Omega_r$ then $\varphi\in D(H_{\!F})$.

Thirdly, if $\varphi\in D(H_{\!F})\subset D(h) $ with $\supp\varphi\subset \Omega_r$ and $\psi\in C_c^\infty(\Omega)$  with $\psi\varphi=0$ 
then $(H_{\!F}\varphi,\psi)=h(\varphi, \psi)=0$ by locality. 
Therefore $\supp H_{\!F}\varphi\subseteq \supp\varphi$.
Similarly  $\supp (H_{\!u})_{\!F}\varphi\subseteq \supp\varphi$.

Finally if $\psi\in C_c^\infty(\Omega_t)$  then
\[
(\psi, H_{\!F}\varphi)=(H_{\!F}\psi, \varphi)=(H_{\!u}\psi, \varphi)=(\psi,  (H_{\!u})_{\!F}\varphi)
\;.
\]
Consequently $ H_{\!u}\varphi=H_{\!F}\varphi$.
\hfill$\Box$

\bigskip

The point of the introduction of the operator $H_u$ is the following key observation.

\begin{lemma}\label{lloc2.1}
The operator $H_{\!u}=-\divv(C_u\nabla)$ is essentially self-adjoint on $C_c^\infty(\Ri^d)$.
In particular $\varphi\in D(H_{\!u})$ if and only if $\varphi\in D((H_{\!u})_{\!F})$ and then $H_{\!u}\varphi=(H_{\!u})_{\!F}\varphi$.
\end{lemma}
\proof\
It follows by construction that the coefficients of  $C_u$ are Lipschitz continuous functions over $\Ri^d$,  $C_u=C$ on $\Omega_t$
and $C_u\geq \mu_u I$ on $\Omega_t^{\,\rm c}$.
In particular $C_u>0$ on $\Ri^d$.
Then $H_{\!u}$ is essentially self-adjoint by the general results cited in the introduction (see \cite{Dav14}, Section~3).
Alternatively it is  a direct corollary of Proposition~6.1 of \cite{RSi4} or Corollary~4.12 in \cite{ERS5}.\hfill$\Box$

\bigskip

It is now a straightforward corollary of  Proposition~\ref{pun2.11} that $H_{\!F}=H$ on the interior sets $\Omega_r$.

\begin{cor}\label{cun2.1}
If $\varphi\in D(H_{\!F})$ and $\supp\varphi\subset \Omega_r$ for some $r>0$ then $\varphi\in D(H)$
and $H\varphi=H_{\!F}\varphi$.
\end{cor}
\proof\
Fix $u, t$ with  $0<u<t<r$ and let $H_{\!u}$ denote the strongly elliptic operator constructed above.
If $\varphi\in D(H_{\!F})$ and $\supp\varphi\subset \Omega_r$ then $\varphi\in D((H_{\!u})_{\!F})$ and $H_{\!F}\varphi=(H_{\!u})_{\!F}\varphi$ by Statement~\ref{pun2.11-3} in Proposition~\ref{pun2.11}.
But then $\varphi\in D(H_u)$ and $(H_{\!u})_{\!F}\varphi=H_{\!u}\varphi$ by Lemma~\ref{lloc2.1}.
Finally $\varphi\in D(H)$ and $H_{\!u}\varphi=H\varphi$ by Statement~\ref{pun2.11-2}  in  Proposition~\ref{pun2.11}.
Thus, by combination of these observations,  $H_{\!F}\varphi=H\varphi$.
\hfill$\Box$

\subsection{Approximate identities}\label{S2-3}

Corollary~\ref{cun2.1}  establishes that $H_{\!F}\varphi=H\varphi$ for the  $\varphi\in D(H_{\!F})$  supported by the interior sets
$\Omega_r$.
Therefore the proof of self-adjointness, i.e.\ the proof that $H_{\!F}=H$, is reduced to studying the boundary behaviour.
One must verify that $H_{\!F}\varphi=H\varphi$ for those  $\varphi\in D(H_{\!F})$  with $\supp\varphi\subset \Gamma_{\!\!s}$ for some 
$s>0$.
The starting point of  the proof is the construction  of an approximate identity $\{\zeta_n\}$ of positive $W^{2,\infty}$-functions  such that 
$\zeta_n D(H_{\!F})\subseteq D(H)$ and $\|\zeta_n\varphi-\varphi\|_{D(H_{\!F})}\to0$ as $n\to\infty$
with $\|\cdot\|_{D(H_{\!F})}$  the graph norm of $H_{\!F}$.
This is  achieved by constructing an  approximate identity on $L_2(0,\infty)$ and then composing it with the distance function $d_\Gamma$.

An analagous tactic was used in the earlier proofs of $L_1$-uniqueness  but then it sufficed to construct an approximate identity $\{\chi_n\}$ of positive $W^{1,\infty}$-functions  convergent with respect to the norm $\|\cdot\|_{D(h)}$.
Both constructions start from the following observation.

\begin{lemma}\label{lun2.11}
There exist $\chi_n\in W^{1,\infty}(0,\infty)$ such that $0\leq \chi_n\leq 1$, $\chi_n(u)=0$ if $u\in\langle 0, n^{-1}]$, $\chi_n(u)=1$ if $u\geq 1$ and $\chi_n\to1$ pointwise on $\langle 0, \infty\rangle$ as $n\to\infty$.
Moreover,
\[
|\chi_n^{\,\prime}(u)|\leq (\log n)^{-1}\,u^{-1}\]
for all $u\in\langle n^{-1},1\rangle$ and large $n\geq1$.
\end{lemma}
\proof\
The  $\chi_n$ are defined by $\chi_n(u)=0$ if $ u\leq n^{-1}$, $\chi_n(u)=- \log nu/\log n$ if $u\in\langle n^{-1},1\rangle $
and $\chi_n(u)=1$ if $u>1$
where $n\in\Ni$.
\hfill$\Box$

\bigskip

The $\chi_n$ are now used to construct the $\zeta_n$ by variation of an argument given in \cite{RSi3}, Section~4.

\begin{prop}\label{pun2.2}
There exist $\zeta_n\in W^{2,\infty}(0,\infty)$ such that $0\leq \zeta_n\leq 1$, $\zeta_n(u)=0$ if $u\in\langle 0, n^{-1}]$, $\zeta_n(u)=1$ if $u\geq 1$ and $\zeta_n\to1$ pointwise on $\langle 0, \infty\rangle$ as $n\to\infty$.
Moreover, there is an $\alpha>0$ such that 
\[
|\zeta_n^{\,\prime}(u)|\leq \alpha\,(\log n)^{-1}\,u^{-1}\;\;\;\;\; \mbox{ and }\;\;\;\;\; |\zeta_n^{\,\prime\prime}(u)|\leq \alpha\,(\log n)^{-1}\,u^{-2}
\]
for all $u\in\langle n^{-1},1\rangle$ and all large $n$.
\end{prop}
\proof\
Set  $\rho_n=\chi_n^{\,2}$.
The $\rho_n$ are positive $W^{2,\infty}$-functions which converge pointwise to the identity on $\langle0,\infty\rangle$.
Moreover,
$\rho_n^{\,\prime}(u)=2\,u^{-1}\,(\log nu)\,(\log n)^{-2}$
for $u\in [n^{-1},1]$ and $\rho_n^{\,\prime}=0$ on the complement of this interval.
The derivative $\rho_n^{\,\prime}$ is continuous at $n^{-1}$ but $\rho_n^{\,\prime}(1)=2\,(\log n)^{-1}$.
Therefore we introduce $\sigma_n$  by
\[
\sigma_n(u)
= \left\{ \begin{array}{ll}
  0 &\hspace{1cm} \mbox{if } u\in[0,n^{-1}\rangle \; ,  \\[5pt]
\rho_n^{\,\prime}(u)-\rho_n^{\,\prime}(1)(u-n^{-1})/(1-n^{-1})&\hspace{1cm} \mbox{if } u\in[n^{-1}, 1] \; , \\[5pt]
  0 & \hspace{1cm}\mbox{if } u \geq 1\; .
         \end{array} \right.
\]
The modified derivative $\sigma_n$ is now continuous both at $n^{-1}$ and at $1$. 
Finally we define $\zeta_n$~by
\[
\zeta_n(u)=N_n^{-1}\int^u_0 dt\,\sigma_n(t)
\]
with $N_n=\int^1_0dt\, \sigma_n(t)$.
It follows that $\zeta_n\in W^{2,\infty}(0,\infty)$, $\zeta_n(u)=0$  if $u\in[0,n^{-1}]$ and $\zeta_n(u)=1$ if $u\geq 1$.
But $N_n\in[1-2/\!\log n,1]$ and 
\[
\zeta_n^{\,\prime}(u)\geq N_n^{-1}\sigma(u)\geq \rho_n^{\,\prime}(u)-\rho_n^{\,\prime}(1)(u-n^{-1})/(1-n^{-1})
\]
for all $u\in[n^{-1},1]$ and all $n>1$.
Then a straightforward calculation establishes that  $\zeta_n^{\,\prime}\geq 0$ on the interval $[n^{-1},1]$.
Hence the $\zeta_n$ are positive.
Moreover, $\zeta_n\to1$ pointwise on $\langle 0, \infty\rangle$ as $n\to\infty$.

Finally the derivatives of the $\zeta_n$ are zero if $u<n^{-1}$ or $u>1$ but  if $u\in [n^{-1},1]$ then
\begin{equation}
\zeta^{\,\prime}_n(u)=-N_n^{-1}\Big(2u^{-1}\log nu/(\log n)^2-2(u-n^{-1})(1-n^{-1})^{-1}/\log n\Big)
\label{eun2.8}
\end{equation}
and 
\begin{equation}
\zeta_n^{\,\prime\prime}(u)=N_n^{-1}\Big(2u^{-2}(1-\log nu)/(\log n)^2+2(1-n^{-1})^{-1}/\log n\Big)
\;.
\label{eun2.9}
\end{equation}
Since $N_n^{-1}\leq 1-2/\!\log n$ the bounds stated in the proposition follow for all $n\geq3$.
\hfill$\Box$

\bigskip

One can now use the $\zeta_n$ to construct an approximate identity with respect to the graph norm of $H_{\!F}$
by setting $\zeta_{r,n}=\zeta_n\circ(r^{-1}d_\Gamma)$.
But to verify the desired convergence properties it is then necessary to have information on the first and second derivatives 
of $d_\Gamma$.

\subsection{Boundary layer estimates}\label{S2-5}
\noindent In order to follow the approximation procedure outlined above it is necessary to estimate  expressions such as $H_{\!F}(\eta_{r,n}\varphi)$ with $\varphi\in D(H_{\!F})$ and $\eta_{r,n }=\eta_n\circ(r^{-1} d_\Gamma)$ where $\eta_n=\one-\zeta_n$.
One can in principle  use (\ref{eun2.3}), with $\chi$ replaced by $\eta_{r,n }$, to calculate $H_{\!F}(\eta_{r,n }\varphi)$
but this immediately raises a problem.
For example one needs to estimate  $H\eta_{r,n }$.
Formally this is given by 
\begin{equation}
H\eta_{r,n }=
-r^{-2}\,\eta_{r,n }^{\,\prime\prime}\,(\nabla\! d_\Gamma, C\nabla\! d_\Gamma)
-r^{-1}\,\eta_{r,n }^{\,\prime}\,\Tr(\,C(D^2d_\Gamma)\,)
-r^{-1}\,\eta_{r,n }^{\,\prime}\,(\divv C).(\nabla\! d_\Gamma)
\label{eun2.5.00}
\end{equation}
where $(\,\cdot\,,\,\cdot\,)$ and $\Tr(\,\cdot\,)$ denote the scalar product and trace on $l_2(\{1,\ldots, d\})$ and
the matrix $(D^2d_\Gamma)=(\,\partial_k\partial_l d_\Gamma\,)$ is the Hessian  of the distance function.
The difficulty is that $d_\Gamma$ is not generally twice-differentiable. 
It is, however, only necessary to have differentiability in a thin boundary layer $\Gamma_{\!\!r}$
and this does follow if  $\Gamma$ is suitably smooth or regular.
Some basic differentiability properties  follow from  geometric properties of $\Omega$ such as convexity and further  estimates on the second derivatives follow from smoothness of the  boundary.
We briefly summarize the properties that suffice for the subsequent discussion.

\smallskip

First the distance $d_\Gamma$ is a Lipschitz function for all domains $\Omega$ and  $|\nabla\! d_\Gamma|\leq1$.
But $d_\Gamma$ is differentiable at  $x\in \Omega$ if and only if $x$  has a unique nearest point in $n(x)\in \Gamma$.
In this  case $(\nabla\! d_\Gamma)(x)=(x-n(x))/|x-n(x)|$
and $|(\nabla\! d_\Gamma)(x)|=1$.
It follows, however,  from Motzkin's theorem (see, for example, \cite{Hor8},  Theorem~2.1.20, or \cite{BEL}, Theorem~2.2.9) that each  $ x\in \Omega$ has  a unique nearest point in $\Gamma$ if and only if $\Omega=\Ri^d\backslash K$ with $K$ a closed convex set.  
Nevertheless this special form of $\Omega$ is not necessary for all  points in a boundary layer $\Gamma_{\!\!r}$ to have a unique nearest point. 
It  follows, however, in the special case  $\Omega=\Ri^d\backslash K$ with $K$ closed and convex that $d_\Gamma$ is convex on (all open convex subsets of) $\Omega$.
Hence  the Hessian $(D^2d_\Gamma)$ is a positive Radon measure  (see, for example, \cite{EvG} Chapter~6).
Alternatively, if $\Omega$ is convex then $d_\Gamma$ is concave and the Hessian is a negative Radon measure.
Typically the measures contain the Hausdorff measure on the boundary.
Specific examples are given in \cite{DeZ}, Sections~6.5.1 and 7.5.1.

Secondly, it follows from  Gilbarg and Trudinger  \cite{GT}, Appendix~14.6, that if $\Omega$ is a $C^2$-domain  then there is a small positive $u$ such that each 
$x\in \Gamma_{\!\!u}$ has a unique near point $n(x)\in \Gamma$, $d_\Gamma\in C^2(\Gamma_{\!\!u})$  and for each $y\in \Gamma$ there exists a ball $B$ such that
$B\cap\Omega^{\rm c}=y$ with the radii of the balls bounded below by $u$.
The last condition is commonly referred to as the uniform sphere condition.
It implies that the principal curvatures $\kappa_1(y),\ldots, \kappa_{d-1}(y)$  at $y\in \Gamma$  are bounded uniformly by $u^{-1}$.
Therefore the Hessian  defined by the curvatures is uniformly bounded.
Although the discussion in \cite{GT}  is restricted to bounded $C^2$-domains the definitions and principal conclusions can be extended in a
uniform manner to unbounded domains (see \cite{For}, Appendix~B).
The principal curvatures and the Hessian remain uniformly bounded.
Moreover, the eigenvalues of the  Hessian   $(D^2d_\Gamma)$    can be calculated  within $\Gamma_{\!\!u}$.
If $x\in \Gamma_{\!\!u}$ then $((D^2d_\Gamma)(x))$ is unitarily equivalent to the diagonal matrix with eigenvalues
$-\kappa_1(y)(1-d_\Gamma(x)\kappa_1(y))^{-1},\ldots,-\kappa_{d-1}(y)(1-d_\Gamma(x)\kappa_{d-1}(y))^{-1},0$ where $y=n(x)$.
  Details of these properties are given in \cite{GT}.
 Therefore
 \[
\Tr((D^2d_\Gamma)(x))=-\sum^{d-1}_{j=1}\kappa_j(y)(1-d_\Gamma(x)\kappa_j(y))^{-1}
 \]
 for all $x\in\Gamma_{\!\!u}$.
 Since  $\kappa_j(y)<u^{-1}$ for each $j\in\{1,\ldots, d-1\}$ and $y\in\Gamma$ it follows that 
  $1-d_\Gamma(x)\kappa_j(y)>0$ for all $x\in\Gamma_{\!\!u}$ and if $r<u$ then 
 $1-d_\Gamma(x)\kappa_j(y)>1-ru^{-1}$ for all $x\in\Gamma_{\!\!r}$.
 Therefore
   \[
 \Tr((d_\Gamma D^2d_\Gamma)(x))=-d_\Gamma(x)M(y)-d_\Gamma(x)^2\sum^{d-1}_{j=1}\kappa_j(y)^2(1-d_\Gamma(x)\kappa_j(y))^{-1}
  \]
  for all $x\in\Gamma_{\!\!u}$ with $M=\sum^{d-1}_{j=1}\kappa_j$, the unnormalized mean curvature.
The eigenvalues of the Hessian matrix $(|D^2d_\Gamma|) $  are then bounded in the  diagonal representation by the substitutions $\kappa_j\to|\kappa_j|$.
 Therefore one obtains estimates
   \[
  |\Tr((d_\Gamma  D^2d_\Gamma)(x))|\leq    \Tr((d_\Gamma |D^2d_\Gamma|)(x))\leq d_\Gamma(x)\,|M|+2\,d_\Gamma(x)^2\,|M|^2
   \]
   for all $x\in\Gamma_{\!\!r}$ with $r\leq u/2$ where $|M|=\sum^{d-1}_{j=1}\|\kappa_j\|_\infty$.
   Consequently, there is a $\gamma>0$ such that
  \begin{equation}
|(d_\Gamma\nabla^2d_\Gamma)(x)|\leq  \Tr((d_\Gamma |D^2d_\Gamma|)(x))\leq \gamma d_\Gamma(x)\leq \gamma r
  \label{eun2.300}
  \end{equation}
  for all $x\in\Gamma_{\!\!r}$.
 Note that a  similar estimate occurs as Condition~(R) in Section~2 of \cite{FMT1} in the analysis of critical  Hardy--Sobolev  inequalities.

Thirdly, one can   estimate  $d_\Gamma\nabla^2d_\Gamma$ for domains with boundaries of lower dimension.
  The classic example is $\Omega=\Ri^d\backslash\{0\}$.
  Then $\Gamma=\{0\}$ and $\Gamma_{\!\!r}$ is a ball of radius $r$ centred at the origin.
  In addition $d_\Gamma(x)=|x|$.
  Therefore $d_\Gamma(D^2d_\Gamma)=P$ where $P_{kl}(x)=(\,\delta_{kl}-e_ke_l\,)$  with $e_k=x_k/|x|$ and $e_l=x_l/|x|$.
  Note that $P=(\,P_{kl}\,)$ is an orthogonal projection on $\Ri^d$  with trace $d-1$.
  Hence (\ref{eun2.300}) is now replaced  by 
  \begin{equation}
 (d_\Gamma\nabla^2d_\Gamma)(x)= \Tr((d_\Gamma D^2d_\Gamma)(x))= \Tr((d_\Gamma |D^2d_\Gamma|)(x))=(d-1)
   \label{eun2.31}
   \end{equation}
  for all $x\in \Gamma_{\!\!r}$ and in fact for all $x\in\Omega$.

Finally  we consider a family of domains which have a structure intermediate between these two cases.
The following discussion has a number of features in common with the papers \cite{BFT} \cite{FMT} and \cite{FMT1}.
As orientation note that   if $\Omega$ is a $C^2$-domain  then the exterior set  $\Epsilon=(\overline\Omega)^{\rm c}$ is also a $C^2$-domain with the same boundary and $\Omega=\Ri^d\backslash \Epsilon^{\,\rm c}$.
First we fix a factorization $\Ri^d=\Ri^s\times \Ri^{d-s}$ with  $s\in\{1,\ldots, d-1\}$ and a $C^2$-domain $\Pi$ in the subspace $\Ri^s$.
Then we set $\Omega=\Ri^d\backslash\overline \Pi$.
It follows that  $\Gamma=\overline \Pi$ and $s$ is equal to the Hausdorff dimension $d_H$ of $\Gamma$.
Now  one can  choose coordinates  $x=(y,z)\in \Ri^d$ with $y\in \Ri^s$ and $z\in \Ri^{d-s}$ such that 
 $d_\Gamma(x)=(d_{\Pi}(y)^{2}+d_A(z)^2)^{1/2}$ where $d_{\Pi}$ is the Euclidean distance to $\Pi$ in $\Ri^s$
 and $d_A$ is the distance to the affine closure $A=\Ri^s$ of $\Pi$.
 In particular $d_A(z)=|z|$  and $d_\Pi(y)=0$ if $y\in\Pi$.
It follows  from the preceding observation that the Hessian  $(D_{\!y}^2d_{\Pi})$  corresponding to $\Pi$ satisfies analogous properties to those established  in the  first case above within a boundary layer $\Pi_{u}=\{y\in \Pi^{\rm c}: d_\Pi(y)<u\}$.
Now $(\,D_x^2d_\Gamma^{\,2}\,)= (\,D_{\!y}^2d_\Pi^{\,2}\,) +(\,D_{\!z}^2d_A^{\,2}\,)$ and 
$(\,D_x^2d_\Gamma^{\,2}\,)=2\,d_\Gamma (\,D_x^2d_\Gamma\,)+2\,(\,D_xd_\Gamma\times D_xd_\Gamma\,)$ etc.
Therefore
\begin{eqnarray}
d_\Gamma (\,D_x^2d_\Gamma\,)=d_\Pi(\, D_y^2d_\Pi\,)+d_A (\,D_z^2d_A\,)-R
\label{eun2.85}
\end{eqnarray}
where $R=(\,D_xd_\Gamma\times D_xd_\Gamma\,)-(\,D_yd_\Pi\times D_yd_\Pi\,)-(\,D_zd_A\times D_zd_A\,)$.
The Hessian of $d_\Gamma$ is a $d\times d$-matrix on $\Ri^d$ and the Hessians of $d_\Pi$ and $d_A$ are matrices on $\Ri^s$ and $\Ri^{d-s}$, respectively.
But $\Tr(R)=|(\nabla_{\!\!x} d_\Gamma)|^2-|(\nabla_{\!\!y} d_\Pi)|^2-|(\nabla_{\!\!z}d_A)|^2=0$ and $d_A(\,D_z^2d_A\,)=P_z$
where $P_z$ is the projection  $P$ of the previous paragraph but on $\Ri^{d-s}$. 
Thus $\Tr(P_z)=d-s-1$.
Therefore 
\[
\Tr((d_\Gamma D^2d_\Gamma)(x))=\Tr((d_\Pi D_{\!y}^2d_\Pi)(y))+(d-s-1)
\;.
\]
A similar identity  is given in Example~4 of \cite{BFT}.
Hence,  by applying (\ref{eun2.300}) to the right hand side with $\Gamma$ replaced by $\overline \Pi$ and $d_\Gamma(x)$ replaced by $d_\Pi(y)$, one obtains
\begin{equation}
|\Tr((d_\Gamma D^2d_\Gamma)(x))-(d-s-1)|\leq \Tr((d_\Pi |D_{\!y}^2d_\Pi|)(y))\leq \gamma d_\Pi(y)\leq \gamma r
\label{eun2.311}
\end{equation}
for all $x\in \Gamma_{\!\!r}$.

All of the foregoing calculations and estimations are intended in the weak sense.
The necessary differentiability is assured by the $C^2$-assumptions and the product structure.
The various bounds are then valid on annular layers $\Gamma_{\!\!r}\backslash\Gamma_{\!\!s}$
where $0<s<r$.
This suffices for the subsequent purposes and avoids problems with measures on the boundaries.

The three inequalities (\ref{eun2.300}), (\ref{eun2.31}) and (\ref{eun2.311}) appear to be rather different but they can be written in a unified form.
There is, however, a distinction between the first two cases and the third.
In all three cases one has 
\begin{equation}
|\Tr((d_\Gamma D^2d_\Gamma)(x))-(d-d_H-1)|\leq \gamma d_\Gamma(x)\leq \gamma r
\label{eun2.321}
\end{equation}
for all $x\in \Gamma_{\!\!r}$ with no further assumptions on $\Omega$.
For example, if $\Omega$ is a $C^2$-domain then $d_H=d-1$ and (\ref{eun2.321}) coincides with (\ref{eun2.300}).
Alternatively, if $\Omega=\Ri^d\backslash\{0\}$ then $\Gamma=\{0\}$,  $d_H=0$ and 
 (\ref{eun2.321}) coincides with (\ref{eun2.31}) but one can choose $\gamma=0$.
Finally, in the intermediate case $\Gamma=\overline \Pi$ and $d_H=s$. 
Therefore  (\ref{eun2.311}) immediately gives  (\ref{eun2.321}).
Bounds of this form  were introduced as Condition~(R) in Sections~2 and~4 of \cite{FMT1} and were used to establish generalized
Hardy--Sobolev inequalities on $L_p$-spaces.
But in the first and second cases the bounds also follow with  the Hessian $(D^2d_\Gamma)$ replaced by its modulus $(|D^2d_\Gamma|)$.
This stronger form also follows in the intermediate case if the domain $\Pi$ is convex.

\begin{prop}\label{npun2.1}
Assume that {\rm either} $\Omega$ is a $C^2$-domain in $\Ri^d$, {\rm or} $\Omega=\Ri^d\backslash \{0\}$ {\rm or} $\Omega=\Ri^d\backslash \overline \Pi$ where $\Ri^d=\Ri^s\times \Ri^{d-s}$ and $\Pi$ is a convex $C^2$-domain in $\Ri^s$.

In each case there is an $r\in\langle0,1]$ and a $\gamma\geq0$ such that 
\begin{equation}
|\Tr((d_\Gamma |D^2d_\Gamma|)(x))-(d-d_H-1)|\leq \gamma d_\Gamma(x)\leq \gamma r
\label{eun2.322}
\end{equation}
for all $x\in \Gamma_{\!\!r}$. 
\end{prop}
\proof\
In the first case $d-d_H-1=0$ and the bounds follow from (\ref{eun2.300}) with $\gamma>0$.
In the second case $d_H=0$ and the bounds follow from (\ref{eun2.31})  with $\gamma=0$.
Finally in the third case, the Hessian  $(D_x^2d_\Gamma)$ is  a positive-definite $d\times d$-matrix 
by the convexity of $\Pi$ which implies the convexity of $d_\Gamma$.
Moreover, $(D_y^2d_\Pi)$ is  a positive-definite $s\times s$-matrix.
Therefore it follows by taking the trace of (\ref{eun2.85}) that 
\[
(d-d_H-1)\leq \Tr((d_\Gamma D^2d_\Gamma)(x))
= \Tr((d_\Gamma |D^2d_\Gamma|)(x))\leq \ Tr((d_\Pi D_y^2d_\Pi)(y))+(d-d_H-1)
\;.
\]
Hence 
\[
0\leq \Tr((d_\Gamma |D^2d_\Gamma|)(x))-(d-d_H-1)\leq \gamma d_\Gamma(x)\leq \gamma r
\]
for all $x\in \Gamma_{\!\!r}$ by application of (\ref{eun2.300}) to the right hand side with $\Gamma$ replaced by $\overline \Pi$.
\hfill$\Box$

\bigskip

Finally it follows from these estimates that one can bound all the terms occurring in the expression (\ref{eun2.5.00}) for $H\eta_{r,n }$.
For orientation consider the second term, the troublesome term with the Hessian of $d_\Gamma$.
First one has $\eta_{r,n}^{\,\prime}=\eta_n^{\,\prime}\circ (r^{-1}d_\Gamma)$ and it follows from Proposition~\ref{pun2.2}
that one has a bound $r^{-1}\,\eta_{r,n }^{\,\prime}(x)\leq \alpha\,(\log n)^{-1}d_\Gamma(x)^{-1}$ if  $x\in\Gamma_{\!\!r}$.
Secondly, it follows from  (\ref{euns2.2}) that $\|C(x)\|\leq \tau_r\,a(x)d_\Gamma(x)^{\delta}$ on $\Gamma_{\!\!r}$.
But $|\Tr(CD^2d_\Gamma)|\leq \|C\|\,\Tr(|D^2d_\Gamma|)$.
Consequently, appealing to (\ref{eun2.322}), one obtains the estimate
\[
|r^{-1}\,\eta_{r,n }^{\,\prime}(x)\,\Tr(\,C(x)(D^2d_\Gamma)(x)\,)|\leq \alpha(\log n)^{-1}\,\tau_r\,(d-d_H-1+\gamma r)\,a(x)d_\Gamma(x)^{\delta-2}
\]for all $x\in\Gamma_{\!\!r}$.
It is also straightforward to derive similar bounds for the first and third terms on the right hand side of (\ref{eun2.5.00}).
For example,
\[
|r^{-2}\,\eta_{r,n }^{\,\prime\prime}(x)\,(\nabla\! d_\Gamma, C\nabla\! d_\Gamma)(x)|
\leq \alpha(\log n)^{-1}\,\tau_r\,\,a(x)d_\Gamma(x)^{\delta-2}
\]
and 
\[
|r^{-1}\,\eta_{r,n }^{\,\prime}(x)\,((\divv C).(\nabla\! d_\Gamma))(x)|
\leq  \alpha(\log n)^{-1}\,\rho_r\,u^{-1}\,a(x)d_\Gamma(x)^{\delta-2}
\;.
\]
Therefore one concludes that 
\begin{equation}
\|(H\eta_{r,n})\varphi\|_2\leq c\,(\log n)^{-1}\,\|a\,d_\Gamma^{\,\delta-2}\varphi\|_2
\label{eun2.33}
\end{equation}
for all $\varphi\in L_2(\Omega)$.

\begin{remarkn}\label{run2.1}
In the intermediate case $\Omega=\Ri^d\backslash\overline \Pi$ convexity was used to establish that the Hessian $(D^2d_\Gamma)$ is positive.
But if $C$ is  a multiple of the identity, e.g.\ if $C=ad_\Gamma^\delta I$, then the positivity is not necessary in bounding $H\eta_{r,n}$.
One may use the estimate (\ref{eun2.321}) in place of (\ref{eun2.322}).
Convexity of $\Pi$ is irrelevant.
\end{remarkn}

It is evident from (\ref{eun2.33}) that if $\delta\geq 2$ then $H\eta_{r,n}$ is a bounded multiplication operator.
If, however,  $\delta<2$ then a problem remains.
It is necessary to establish that if $\varphi\in D(H_{\!F})$ and $\supp\varphi\subseteq \Gamma_{\!\!r}$ then $\varphi\in 
D(d_\Gamma^{\,\delta-2})$.
This requires a form of the Rellich inequality.

\subsection{Hardy--Rellich boundary estimates}\label{S2-4}

In an earlier paper \cite{Rob12} we established Rellich inequalities for the Friedrichs' extension for a general class
of degenerate elliptic operators.
The derivation, which was based on ideas of Agmon \cite{Agm1} and Grillo \cite{Gril}, started from a Hardy inequality and depended
on estimates on $\nabla^2d_\Gamma$.
But to derive the Rellich inequalities on a boundary layer $\Gamma_{\!\!r}$ one  only needs the Hardy inequality and the distance estimates 
on $\Gamma_{\!\!r}$.
Since we have already obtained estimates on $\nabla^2d_\Gamma$ it is now necessary to establish Hardy inequalities under matching assumptions.
We begin by introducing comparison operators and  forms.

The matrix $C$ of coefficients of $H$ is equivalent to the diagonal matrix $a d_\Gamma^{\,\delta}I$
on the boundary layer $\Gamma_{\!\!r}$ by (\ref{euns2.2}).
Therefore we  introduce the  operator $H_\delta=-\sum^d_{k=1}\partial_k\,c\,\partial_k$ 
with $c$ a strictly, positive, bounded Lipschitz continuous function.  
Moreover, $c= a d_\Gamma^{\,\delta}$ on $\Gamma_{\!\!r}$.
The behaviour of $c$ in the interior is irrelevant to the following boundary layer estimates.
We denote the corresponding form by $h_\delta$. 
Both the operator and form are closable and for simplicity  we  use $H_\delta$ and $h_\delta$ to denote the corresponding closures.
Now we derive a  Hardy inequality for $h_\delta$ on  functions supported in $\Gamma_{\!\!r}$.
Similar boundary estimates are given in \cite{FMT1}.

\begin{prop}\label{pun3.0}
Assume  that either $\Omega$ is  a $C^2$-domain in $\Ri^d$, or $\Omega=\Ri^d\backslash\{0\}$, or $\Omega=\Ri^d\backslash\overline \Pi$
with $\Pi$ a  $C^2$-domain in the subspace $\Ri^s$.
It follows that   there is a $c_a>0$ such that if $\alpha_{1,r}=d-d_H+\delta-2-c_ar>0$, with $d_H$ the Hausdorff dimension of $\Gamma$,
then
\begin{equation}
h_\delta(\varphi)\geq (\alpha_{1,r}/2)^2\,\|a^{1/2}d_\Gamma^{\,\delta/2-1}\varphi\|_2^2 
\label{eun2.10}
\end{equation}
for all $\varphi\in D(h_\delta)$ with $\supp\varphi\subset \Gamma_{\!\!r}$ and for all $\delta\geq0$ and $r\in\langle0,1]$.
\end{prop}
\proof\
The proof follows standard reasoning. 
Set $\chi=d_\Gamma^{\,\delta-1}\nabla \!d_\Gamma$.
Then
\[
\divv\chi=(\delta-1+d_\Gamma\nabla^2d_\Gamma) \,d_\Gamma^{\,\delta-2}
\geq (d-d_H+\delta-2-\gamma r)\,d_\Gamma^{\,\delta-2}
\]
on $\Gamma_{\!\!r}$ where we have used (\ref{eun2.321}). 
Convexity of $\Pi$ is not necessary for this estimate.
Thus if $t<r$ and $\varphi\in C_c^\infty(\Gamma_{\!\!t})$ then
\[
(d-d_H+\delta-2-\gamma r)\,\int_{\Gamma_{\!\!r}}a d_\Gamma^{\,\delta-2}\varphi^2
\leq \int_{\Gamma_{\!\!r}}a\,(\divv\chi)\varphi^2
\leq \int_{\Gamma_{\!\!r}}|\chi.\nabla a|\,\varphi^2+2\int_{\Gamma_{\!\!r}}a\,|\chi.\nabla\varphi|\,|\varphi|
\]
by partial integration.
But $|\chi.\nabla a|\leq (\gamma_a r)\,a d_\Gamma^{\,\delta-2}$ on $\Gamma_{\!\!r}$ where $\gamma_a$ is the $L_\infty$-norm of 
$ |\nabla a|/a$ on $\Gamma_{\!\!r}$.
Therefore
\[
(d-d_H+\delta-2-c_a r)\,\int_{\Gamma_{\!\!r}}a\,d_\Gamma^{\,\delta-2}\varphi^2
\leq 2\int_{\Gamma_{\!\!r}}a\,|\chi.\nabla\varphi|\,|\varphi|
\]
with $c_a=\gamma+\gamma_a$.
Then setting $\psi=d_{\Gamma}^{\,-\delta/2+1}$ one has
\[
\Big(\int_{\Gamma_{\!\!r}}a\,|\chi.\nabla\varphi|\,|\varphi|\Big)^2
\leq \int_{\Gamma_{\!\!r}}ad_\Gamma^{\,\delta}|\nabla\varphi|^2\int_{\Gamma_{\!\!r}}ad_\Gamma^{\,\delta-2}\varphi^2
\;.
\]
Combination of the  last two estimates gives (\ref{eun2.10}) for all  $\varphi\in C_c^\infty(\Gamma_{\!\!t})$  and then by closure
for all  $\varphi\in D(h_\delta)$ 
\hfill$\Box$ 

 \bigskip

Next we derive a boundary Rellich inequality from the Hardy inequality (\ref{eun2.10}) essentially as a corollary of the arguments 
of \cite{Rob12}.
Some modifications of the previous arguments have to be made since we are only dealing with inequalities close to the boundary.
First we remark that the bounds (\ref{euns2.2}) and (\ref{euns2.2.1}) lead to the equivalence
\[
\sigma_{\!r}\,h_\delta(\varphi)\leq h(\varphi)\leq \tau_{\!r}\,h_\delta(\varphi)
\]
for all $\varphi\in D(h)=D(h_\delta)$ with $\supp\varphi\subseteq \Gamma_{\!\!r}$.
Therefore the Hardy inequality (\ref{eun2.10})  for $h_\delta$ gives a similar inequality  for the form $h$,  if $\alpha_{1,r}=d-d_H+\delta-2-c_ar>0$, 
but the constant $(\alpha_{1,r}/2)^2$ is now multiplied by a factor $\sigma_{\!r}/\tau_{\!r}$.
It is convenient to set $\upsilon_r=\tau_{\!r}/\sigma_{\!r}$.
Note that $\sigma_{\!r}, \tau_{\!r}, \upsilon_r\to1$  as $r\to0$.

\begin{prop}\label{pex2.1}
Adopt the assumptions of Proposition~$\ref{pun3.0}$ with $\delta\in[0,2\rangle$ and $r\in\langle0,1]$.
Set $\alpha_{1,r}=d-d_H+\delta-2-c_ar>0$.

If  $\alpha_{1,r}^2>\upsilon_{\!r}(2-\delta+\gamma_a r)^2$ where $\gamma_a$ is the $L_\infty(\Gamma_{\!\!r})$-norm of $|\nabla a|/a$
then
\[
\|H_{\!F}\varphi\|_2^2\geq (\alpha_{2,r}/4)^2\|a\,d_\Gamma^{\,\delta-2}\varphi\|_2^2
\]
for all $\varphi\in D(H_{\!F})$ with $\supp\varphi\subset\Gamma_{\!\!r}$
where  $\alpha_{2,r}=\sigma_{\!r}(\alpha_{1,r}^2-\upsilon_{\!r}(2-\delta+\gamma_a r)^2)$.
\end{prop}

\begin{remarkn}\label{run2.2} The  condition $\alpha_{1,r}^2>\upsilon_{\!r}(2-\delta+\gamma_a r)^2$  of the proposition can be clarified by noting that in the limit $r\to0$ it is equivalent to
$(d-d_H+\delta-2)^2>(2-\delta)^2$.
Since $d\geq d_H$ and $2-\delta>0$ this is equivalent to  $d-d_H+2\delta-4>0$.
But this is  the strict form of the condition (\ref{eun1.3}) under which we aspire to prove self-adjointness of $H$.
This indicates  the linkage between existence of  the Rellich inequality and self-adjointness.
Moreover, in this limit the Rellich constant assumes the familiar form $(d-d_H)^2(d-d_H+2\delta-4)^2/16$.
\end{remarkn}

\begin{remarkn}\label{run2.11} The  condition $\alpha_{1,r}^2>\upsilon_{\!r}(2-\delta+\gamma_a r)^2$  can also be replaced by an explicit condition on $\delta$.
Since $2-\delta>0$,  $r\leq1$  and $\upsilon_r\to1$ as $r\to0$ 
 it suffices that   
$\delta\geq \delta_0+\delta_1 r$ with $\delta_0=2-(d-d_H)/2$,  $\delta_1>0$ and $r\in\langle0,r_0]$ where $\delta_1r_0\leq 2-\delta_0$.
\end{remarkn}

\medskip
\noindent{\bf Proof of Proposition~\ref{pex2.1}}$\;$ 
The proof is an  elaboration of the argument used to prove 
Theorem~1.2 in \cite{Rob12}.
In particular it uses the approximate identity $\rho_n$ constructed in Proposition~3.1 of that reference.
This corresponds to the definition $\rho_n=\chi_n\circ d_\Gamma$ where $\chi_n$ is given by Lemma~\ref{lun2.11}.
It should not be confused with the $\rho_n$ used in the proof of Proposition~\ref{pun2.2}.

First,  setting $\beta=(\alpha_{1,r}/2)(a^{1/2}d_\Gamma^{\,\delta/2-1})$ and $\beta_m=\beta(1+\beta/m)^{-1}$, with $m\geq 1$,
it follows from (\ref{eun2.10})  that 
\begin{equation}
h_\delta(\varphi)\geq \|\beta\varphi\|_2^2\geq \|\beta_m\varphi\|_2^2
\label{ex2.22}
\end{equation}
for all $\varphi\in D(h_\delta)$.
The $\beta_m$ form a  sequence of  bounded functions on $\Omega$ which converges monotonically upward to $\beta$
as $m\to\infty$.
Further,  set $\beta_{n,m}=\rho_n\beta_m$ with $\rho_n$ the approximate identity of  Proposition~3.1 in \cite{Rob12}.
In particular it follows from the proposition that  $0\leq\rho_n\leq1$,  $\supp\rho_n\subset \Omega_{1/n}$, $\rho_n\to1$ and $(\varphi, \Gamma(\rho_n)\varphi)\to 0$ as $n\to\infty$ for all $\varphi\in D(h)$.

Secondly, fix $\varphi\in D(H_{\!F})$ with  $\supp\varphi\subset \Gamma_{\!\!r}$ where $r<1$.
If  $\varphi_p=\varphi(1+\varphi/p)^{-1}$ with $p\geq1$
then  $\varphi_p\in D(h)\cap L_\infty(\Omega)$, $\supp\varphi_p=\supp\varphi \subset \Gamma_{\!\!r}$ and $\|\varphi_p-\varphi\|_{D(h)}\to0$ as $p\to\infty$, by Lemma~2.6 of \cite{Rob12}.
Moreover, $\supp\beta_{n,m}\varphi_p\subseteq  \Gamma_{\!\!r}\cap\Omega_{1/n}$ and $\beta_{n,m}\varphi_p\in D(h)$.

The starting point of the proof is the Dirichlet form  identity
\[
h(\varphi,\beta_{n,m}^{\,2}\varphi_p)
=h(\beta_{n,m}\varphi_p)-(\varphi_p,\Gamma(\beta_{n,m})\varphi_p)
-h(\varphi-\varphi_p,\beta_{n,m}^{\,2}\varphi_p)
\]
(see \cite{Rob12}) which immediately leads to the estimate
\begin{equation}
(H_{\!F}\varphi,\beta_{n,m}^{\,2}\varphi_p)\geq \sigma_{\!r}\,h_\delta(\beta_{n,m}\varphi_p)-
\tau_{\!r}\,(\varphi_p,a d_\Gamma^{\,\delta}\,|\nabla\!\beta_{n,m}|^2\,\varphi_p)
-h(\varphi-\varphi_p)^{1/2}\,h(\beta_{n,m}^{\,2}\varphi_p)^{1/2}
\;.
\label{ex2.23}
\end{equation}
But $|(H_{\!F}\varphi,\beta_{n,m}^{\,2}\varphi_p)|\leq \|H_{\!F}\varphi\|_2\,\|\beta_{m}\beta\varphi_p\|_2$ since $\rho_n\in\langle0,1]$
and $\beta_m\leq \beta$.
Moreover, 
\[
h_\delta(\beta_{n,m}\varphi_p)\geq \|\beta_{n,m}\,\beta\varphi_p\|_2^2= \|\rho_n\beta_{m}\beta\varphi_p\|_2^2
\]
 by the Hardy inequality (\ref{ex2.22}).
Thus inserting  these estimates into (\ref{ex2.23}) and taking the limit $n\to\infty$ gives
\begin{eqnarray}
\|H_{\!F}\varphi\|_2\,\|\beta_{m}\,\beta\varphi_p\|_2&\geq & \sigma_{\!r}\,\|\beta_{m}\beta\varphi_p\|_2^2
-\tau_{\!r}\limsup_{n\to\infty}(\varphi_p,a d_\Gamma^{\,\delta}\,|\nabla\!\beta_{n,m}|^2\,\varphi_p)\nonumber\\[5pt]
&&\hspace{1cm}
-h(\varphi-\varphi_p)^{1/2}\,\limsup_{n\to\infty}h(\beta_{n,m}^{\,2}\varphi_p)^{1/2}
\,.
\label{ex2.24}
\end{eqnarray}
Now we successively consider the limits  $n\to\infty$, $p\to\infty$ and $m\to\infty$.

\smallskip

The $n$-dependence of second term on the right hand side of (\ref{ex2.24}) can be handled by using the Leibniz rule for 
$\nabla\!\beta_{n,m} =\nabla(\rho_n\beta_m)$ and squaring.
There are three terms.
First there is a  term which depends on $n$ through a factor $\rho_n^{\,2}$ and since $\rho_n\to\one$ in the  limit  $n\to\infty$ it leads to a contribution  $\tau_{\!r}(\varphi_p,a d_\Gamma^{\,\delta}|\nabla\!\beta_{m}|^2\varphi_p)$.
Secondly  there is a  correction term,   which depends on $|\nabla\!\rho_n|^{2}$, given by 
$\tau_{\!r}(\varphi_p,a d_\Gamma^{\,\delta}\,|\nabla\!\rho_n|^{2}\,\beta_m^{\,2}\, \varphi_p)$.
But this is bounded by 
$ \tau_{\!r}\|\beta_m\|_\infty^2\,(\varphi_p,a d_\Gamma^{\,\delta}\,|\nabla\!\rho_n|^{2}\,\varphi_p)$
and  the latter expression tends to zero as $n\to\infty$ since $(\psi,\Gamma(\rho_n)\psi)\to0$ for all $\psi\in D(h)$
 (see \cite{Rob12}, Proposition~3.1).
 Here it is essential that $\varphi_p\in D(h)$.
  Finally the cross-terms in the correction also tend to zero by an application of the Cauchy--Schwarz inequality.

Next consider the limit $n\to\infty$ in the last term on the right hand side of (\ref{ex2.24}).
Since  $h(\beta_{n,m}^{\,2}\varphi_p)=h(\rho_n^{\,2}\beta_m^{\,2}\varphi_p)$
this can again be handled by using  the Leibniz rule  and arguing as above.
One concludes that $\limsup_{n\to\infty}h(\beta_{n,m}^{\,2}\varphi_p)\leq h(\beta_m^{\,2}\varphi_p)$.
But  one can reapply the Leibniz rule to deduce a further $p$-independent  bound. 
There are a variety of terms   but they all have a bound $a_m (h(\varphi_p)+\|\varphi_p\|_2^2)$ with $a_m>0$ independent of~$p$.
For example one has a term with a factor $\beta_m^{\,4}$ and this can be bounded by  $\|\beta_m\|_\infty^{4}\leq m^4$.
Alternatively  there is a term with a factor $\beta_m^{\,2}\,d_\Gamma^{\,\delta}\,|\nabla\!\beta_m|^2$ and this is also bounded by a multiple of $m^4$ by direct calculation of~$\nabla\!\beta_m$.
Finally $\|\varphi_p\|_{D(h)}\leq \|\varphi\|_{D(h)}$ and so one has a  bound uniform in $p$ of the form
$\kappa\,m^4\, \|\varphi\|_{D(h)}^2$.
Therefore the basic inequality becomes
\begin{eqnarray}
\|H_{\!F}\varphi\|_2\,\|\beta_{m}\,\beta\varphi_p\|_2&\geq & \sigma_{\!r}\,\|\beta_m\beta\varphi_p\|_2^2
-\tau_{\!r}(\varphi_p,a d_\Gamma^{\,\delta}\,|\nabla\!\beta_m|^2\,\varphi_p)\nonumber\\[5pt]
&&\hspace{1cm}
-\kappa \,m^4\,h(\varphi-\varphi_p)^{1/2}\, \|\varphi\|_{D(h)}^2
\,.
\label{ex2.242}
\end{eqnarray}
Now we consider the $p\to\infty$ limit.

First the final term on the right  of (\ref{ex2.242})  tends to zero because $\|\varphi_p-\varphi\|_{D(h)}\to 0$ as $p\to\infty$.
Secondly $\|\beta_{m}\,\beta\varphi_p\|_2^2\leq m^2\|\beta\varphi_p\|_2^2\leq m^2\,h_\delta(\varphi_p)$ by the Hardy inequality (\ref{ex2.22}).
Therefore $\|\beta_{m}\,\beta\varphi_p\|_2\to
\|\beta_m\beta\varphi\|_2$ as $p\to\infty$ for the same reason.
It remains to examine the second term on the right.
This is more complicated and requires some calculation.

It follow by definition that $\nabla\!\beta_m=(1+\beta/m)^{-2}\nabla\!\beta$. 
Then since $\beta$ is proportional to the product of  $a^{1/2}$
and $d_\Gamma^{\,\delta/2-1}$ one can again use the Leibniz rule.
The leading term is given by the square of the contribution in which the derivatives are on the $d_\Gamma$.
This can be calculated exactly as in \cite{Rob12} and one obtains an estimate
\[
\sigma_{\!r}(\upsilon_r(2-\delta)^2/\alpha_{1,r}^2)(\varphi_p,\beta_m^4\varphi)
\leq \sigma_{\!r}(\upsilon_r(2-\delta)^2/\alpha_{1,r}^2)\,\|\beta_m\beta\varphi_p\|_2^2
\;.
\]
The first correction is the square of the contribution with the derivatives on $a^{1/2}$.
A straightforward calculation gives an estimate
\[
\sigma_{\!r}(\upsilon_r(2-\delta)^2/\alpha_{1,r}^2)\,\gamma_a^2\,(r/(2-\delta))^2\,\|\beta_m\beta\varphi_p\|_2^2
\]
where $\gamma_a$ is again the $L_\infty$-norm of $|\nabla a|/a$ on $\Gamma_{\!\!r}$.
Then there are the cross terms which can be estimated with the Cauchy--Schwarz inequality by the usual $(\varepsilon, \varepsilon^{-1})$-method.
Choosing $\varepsilon =r\,(\gamma_a/(2-\delta))$ the sum of all the contributions is then given by
\[
\sigma_{\!r}(\upsilon_r(2-\delta)^2/\alpha_{1,r}^2)\,(1+\gamma_a r/(2-\delta))^2\,\|\beta_m\beta\varphi_p\|_2^2
\]
and the final factor again converges to $\|\beta_m\beta\varphi\|_2^2$ as $p\to \infty$.
It should be emphasized  that $\|\beta_m\beta\varphi\|_2^2\leq m^2\,h_\delta(\varphi)\leq (m^2/\sigma_r)\,h_\delta(\varphi)$ for all $\varphi\in D(h)$, 
the equivalence of $h_\delta$ and $h$  and    the Hardy inequality (\ref{ex2.22}).
In particular for $m$ fixed  $\|\beta_m\beta\varphi\|_2<\infty$ for all $\varphi\in D(H_{\!F})\subseteq D(h)$.

Combining all these observations one  deduces that after the $p\to \infty$ limit in (\ref{ex2.242}) one can divide by a common finite  factor
$\|\beta_m\beta\varphi\|_2$ to obtain
\begin{equation}
\|H_{\!F}\varphi\|_2\geq  \sigma_{\!r}\,(1-\upsilon_r(2-\delta+\gamma_a r)^2/\alpha_{1,r}^2)\,
\|\beta_m\beta\varphi\|_2^2
\label{eun2.243}
\end{equation}
for all $\varphi\in D(H_{\!F})$ with support in $\Gamma_{\!\!r}$.
Then the first factor on the right is strictly positive if $\alpha_{1,r}^2>\upsilon_r(2-\delta+\gamma_a r)^2$.
Finally one can take the limit $m\to\infty$.
The limit now exists by the monotone convergence theorem and $\|\beta_m\beta\varphi\|_2^2\to \|\beta^2\varphi\|_2^2$.
Therefore  (\ref{eun2.243}) immediately gives 
 the Rellich inequality  of the proposition with the stated value of $\alpha_{2,r}$.
\hfill$\Box$

\bigskip

The Rellich inequalities of Proposition~\ref{pex2.1} play a crucial role in the subsequent proof of self-adjointness of $H$.
But it is critical that they are valid for all $\varphi\in D(H_{\!F})$ with support in a thin boundary layer $\Gamma_{\!\!r}$.
Most of the literature devoted to the  Rellich inequality is restricted to interior estimates on subspaces of smooth  functions such as  
$C_c^\infty(\Omega)$ or  $ W^{2,2}_0(\Omega)$ (see, for example, \cite{BEL} Chapter~6 and references therein).
But these subspaces are a core of $H_{\!F}$ if and only if  the restriction of $H$ to the subspace is essentially self-adjoint.

\section{Self-adjointness}\label{S3}

After these extensive preparations we turn to the main focus of this article, criteria for  the  essential self-adjointness of the degenerate operators $H=-\divv(C\nabla)$ on $L_2(\Omega)$.
First we establish sufficiency conditions for  $\Omega$  in one of the three classes analyzed in Subsections~\ref{S2-5} and \ref{S2-4}. 
Subsequently we discuss necessary conditions largely for $C^2$-domains.
 
 \begin{thm}\label{tun5.1}
Assume  that either $\Omega$ is  a $C^2$-domain in $\Ri^d$, or $\Omega=\Ri^d\backslash\{0\}$, or $\Omega=\Ri^d\backslash\overline \Pi$
with $\Pi$ a convex $C^2$-domain in the subspace $\Ri^s$.
Further assume that the coefficients $C$ of the elliptic operator $H$ satisfy  conditions $(\ref{euns3.2})$ and $(\ref{eun3.21})$.

If $\delta>2-(d-d_H)/2$ where $d_H$ is the Hausdorff dimension of the boundary $\Gamma$ of $\Omega$ then  $H$ is self-adjoint,
i.e.\ $C_c^\infty(\Omega)$ is a core of $H_{\!F}$.
\end{thm}
\proof\
 Since we assume that $H$ is closed it suffices  for   self-adjointness to establish that $D(H_{\!F})\subseteq D(H)$.
Fix   $\varphi\in D(H_{\!F})$ and $\chi\in W^{2,\infty}(\Omega)$ with $0\leq \chi\leq 1$,  $\chi=1$ on $\Omega_t$ and $\chi=0$ on $\Gamma_{\!\!u}$ where $u<t<r$.
Then $\varphi=(\one-\chi)\varphi+\chi\varphi$ 
 and $\chi\varphi\in D(H)$ by Corollary~\ref{cun2.1}.
 Thus to deduce that $\varphi\in D(H)$ it suffices to prove that $(\one-\chi)\varphi\in D(H)$.
 But $\supp(\one-\chi)\varphi\subseteq \Gamma_{\!\!t}$.
Hence  in the remainder of the proof we may effectively assume $\supp\varphi\subset \Gamma_{\!\!t}$ for some $t\in\langle0,r\rangle$.

Next set $\zeta_{r,n}=\zeta_n\circ(r^{-1}d_\Gamma)$ and $\varphi_{r,n}=\zeta_{r,n}\,\varphi$ where $\zeta_n$ is  the approximate
identity constructed in Proposition~\ref{pun2.2}.
 It follows that $\supp\varphi_{r,n}\subseteq \Gamma_{\!\!r}$.
Now  $\zeta_{r,n}\in W^{2,\infty}(\Omega)$ and $\varphi_{r,n}\in D(H_{\!F})$ by  Proposition~\ref{pun2.1}.
 But it also follows from Corollary~\ref{cun2.1} that $\varphi_{r,n}\in D(H)$.
 Moreover, setting  $\eta_{r,n}=\one-\zeta_{r,n}$, one has  $\eta_{r,n}D(H_{\!F})=(\one-\zeta_{r,n})D(H_{\!F})\subseteq D(H_{\!F})$ and 
 \[
H_{\!F}(\varphi-\varphi_{r,n})=H_{\!F}(\eta_{r,n}\varphi)=\eta_{r,n}(H_{\!F}\varphi)+(H\eta_{r,n})\varphi -2\,
\Gamma(\eta_{r,n},\varphi)\;.
\]
Note that $\eta_{r,n}=1$ if $d_\Gamma(x)\leq r n^{-1}$,  $\eta_{r,n}=0$ if $d_\Gamma(x)\geq r$ and $\eta_{r,n}\to0$ pointwise as $n\to\infty$.
Therefore $\supp \eta_{r,n}$ is not necessarily compact.
Hence it is not evident that $\Gamma(\eta_{r,n},\varphi)$ is square-integrable.
But this follows from the identity since all other terms are in $L_2(\Omega)$ for $r$ and $n$ fixed.
Then, however,
\begin{eqnarray}
\|H_{\!F}(\varphi-\varphi_{r,n})\|_2\leq \|\eta_{r,n}(H_{\!F}\varphi)\|_2+\|(H\eta_{r,n})\varphi\|_2
 +2\,\|\Gamma(\eta_{r,n},\varphi)\|_2
 \;. \label{eun3.1}
 \end{eqnarray}
 Now we use this bound to establish that  $\|H_{\!F}(\varphi-\varphi_{r,n})\|_2\to0$ as $n\to\infty$ through term by term consideration of the right hand side.
 
 The first term on the right of (\ref{eun3.1}) tends to zero as $n\to\infty$ by dominated convergence since $(H_{\!F}\varphi)^2\in L_1(\Omega)$ and  $\eta_{r,n}\to0$ pointwise.  
This argument is independent of   the assumption $\delta>2-(d-d_H)/2$.
 But this condition is of importance in the consideration of the other two terms on the right of (\ref{eun3.1}).
 
The estimation of the second term depends on validity of the Rellich inequality of Proposition~\ref{pex2.1} on $\Gamma_{\!\!r}$
and this requires the condition $\alpha_{1,r}^2>\upsilon_{\!r}(2-\delta+\gamma_a r)^2$.
But if $r$ is sufficiently small this follows from an explicit bound $\delta\geq \delta_0+\delta_1r$ where $\delta_0=2-(d-d_H)/2$ and $\delta_1>0$ by Remark~\ref{run2.11}.
Therefore we start by making this assumption.
Hence one can  conclude from the Rellich inequality  that if  $\varphi\in D(H_{\!F})$ with $\supp\varphi\subseteq \Gamma_{\!\!r}$ then $\varphi\in D(d_\Gamma^{\,\delta-2})$
and $\|a\,d_\Gamma^{\,\delta-2}\varphi\|_2\leq (4/\alpha_{2,r})\|H_{\!F}\varphi\|_2$.
But then the estimates (\ref{eun2.5.00}) and (\ref{eun2.33}) of Subsection~\ref{S2-5}  give bounds
\[
\|(H\eta_{r,n})\varphi\|_2\leq c\,(4/\alpha_{2,r})\,(\log n)^{-1}\,\|H_{\!F}\varphi\|_2
\;
\]
Therefore $\|(H\eta_{r,n})\varphi\|_2\to0$ as $n\to \infty$ under the proviso $\delta\geq \delta_0+\delta_1r$.

There are two important elements in this argument, the Rellich inequality of Propostion~\ref{pex2.1} and  the estimate (\ref{eun2.33}).
If $\Omega=\Ri^d\backslash\overline\Pi$ then convexity of $\Pi$ is needed to derive (\ref{eun2.33})
although it is not necessary for the Rellich inequality.
It is also not necessary if $C$ is a multiple of the identity (see Remark~\ref{run2.11}).

It now remains to consider the third  term on the right hand side of (\ref{eun3.1}).
But this is estimated by
\[
\|\Gamma(\eta_{r,n},\varphi)\|_2^2\leq \int_{\Gamma_{\!r,n}}{\!\!\!}dx\,\Gamma(\eta_{r,n})\,\Gamma(\varphi)
\leq \tau_r^2\int_{\Gamma_{\!r,n}}{\!\!\!}dx\,\Gamma_{\!\!\delta}(\eta_{r,n})\,\Gamma_{\!\!\delta}(\varphi)
\;,
\]
where $\Gamma_{\!\!\delta}$ is  the {\it carr{\`e} du champ} corresponding to $H_\delta$ and 
$\Gamma_{\!\!r,n}=\{x\in\Omega: rn^{-1}\leq d_\Gamma(x)\leq r\}$.
Then
\[
\int_{\Gamma_{\!r,n}}{\!\!\!}dx\,\Gamma_{\!\!\delta}(\eta_{r,n})\,\Gamma_{\!\!\delta}(\varphi)
\leq r^{-2}\int_{\Gamma_{\!r,n}}{\!\!\!}dx\,(a d_\Gamma^{\,\delta}\, \eta_{r,n}^{\,\prime})^2\,|\nabla\varphi|^2
\;.
\]
Now let $J_n$ denote the last integral.
Then
\begin{eqnarray*}
J_n=r^{-2}\int_{\Gamma_{\!r,n}}{\!\!\!}dx\,(a d_\Gamma^{\,\delta}\,  \eta_{r,n}^{\,\prime})^2\,|\nabla\varphi|^2
=r^{-2}\int_{\Gamma_{\!r,n}}{\!\!\!}dx\,(\nabla\varphi).((a d_\Gamma^{\,\delta}\, \eta_{r,n}^{\,\prime})^2\,\nabla\varphi)
\;.
\end{eqnarray*}
Integrating by parts 
one  obtains
\[
J_n=-r^{-2}\int_{\Gamma_{\!r,n}}{\!\!\!}dx\,\varphi\,ad_\Gamma^{\,\delta}(\eta_{r,n}^{\,\prime})^2\,\divv(ad_\Gamma^{\,\delta}\,\nabla\varphi)
-r^{-2}\int_{\Gamma_{\!r,n}}{\!\!\!}dx\,\varphi\,(\nabla (ad_\Gamma^{\,\delta}(\eta_{r,n}^{\,\prime})^2)).(ad_\Gamma^{\,\delta}\, \nabla\varphi)
\;.
\]
Denote the two terms on the right hand side by $J_{1,n}$ and $J_{2,n}$, respectively.
Then it follows from the  Cauchy--Schwarz inequality that 
\[
|J_{1,n}|^2\leq \Big( r^{-4}\int_{\Gamma_{\!r,n}}{\!\!\!}dx\,|(ad_\Gamma^{\,\delta}(\eta_{r,n}^{\,\prime})^2)|^2 \,|\varphi|^2\Big) \,\|H_{\!F}\varphi\|_2^2\;.
\]
Hence  the bounds  on $\eta_n^{\,\prime}$ which follow from Proposition~\ref{pun2.2} give 
\[
|J_{1,n}|^2\leq \alpha^4\,(\log n)^{-4} \Big( \int_{\Gamma_{\!r,n}} dx\,|ad_\Gamma^{\,\delta-2} \varphi|^2 \Big) \,\|H_{\!F}\varphi\|_2^2
\;.
\]
Thus if $\delta\geq2$ then
\[
|J_{1,n}|^2\leq \alpha^4\,(\log n)^{-4}\,\lambda\, \|\varphi\|_2^2 \,\|H_{\!F}\varphi\|_2^2
\]
where $\lambda$ is the upper bound of $a$.
If, however, $\delta\in\langle\delta_0+\delta_1r, 2\rangle$ then 
\begin{eqnarray}
|J_{1,n}|^2\leq \alpha^4(4/\alpha_{2,r})^2\,(\log n)^{-4}\,\| H_{\!F}\varphi\|_2^4
\label{eun3.4}
\end{eqnarray}
where we have  used the Rellich inequality  as above.
Therefore   $J_{1,n}\to 0$ as $n\to\infty$ for all $\delta>\delta_0+\delta_1r$.
 
 Next consider $J_{2,n}$.
One has 
\[
|J_{2,n}|\leq r^{-2}\int_{\Gamma_{\!r,n}} {\!\!}|\varphi|\,|(\nabla (ad_\Gamma^{\,\delta}(\eta_{r,n}^{\,\prime})^2)).(ad_\Gamma^{\,\delta}\nabla\varphi)|
\;.
\]
But 
\begin{eqnarray}
|(\nabla (ad_\Gamma^{\,\delta}(\eta_{r,n}^{\,\prime})^2)).(d_\Gamma^{\,\delta}\nabla\varphi)|
&\leq& a \Big(|\eta_{r,n}^{\,\prime}||\nabla d_\Gamma^{\,\delta}|+2\,d_\Gamma^{\,\delta}\,|\nabla\eta_{r,n}^{\,\prime}|\Big)\,|(ad_\Gamma^{\,\delta}\,\eta_{r,n}^{\,\prime}\nabla\varphi)|\nonumber\\[5pt]
&&{}+2\,\Big((|\nabla a|/a)\,|ad_\Gamma^{\,\delta}\eta_{r,n}^{\,\prime})|\Big)\,|ad_\Gamma^{\,\delta}\eta_{r,n}^{\,\prime}\nabla\varphi|
\label{eun3.41}
\end{eqnarray}
where we have moved one of the  $\eta_{r,n}^{\,\prime}$ from the first factor in this expression to the second.
Now
$|\nabla d_\Gamma^{\,\delta}|\leq \delta\,d_\Gamma^{\,\delta-1}$ since $|\nabla d_\Gamma|\leq1$.
Hence, by another application of the bounds of Proposition~\ref{pun2.2}, one deduces that there is  an $\tilde \alpha>0$ such that 
\[
|(\nabla (a\,d_\Gamma^{\,\delta}(\eta_{r,n}^{\,\prime})^2)).(ad_\Gamma^{\,\delta}\,\nabla\varphi)|
\leq\tilde \alpha\,r(\log n)^{-1} \,d_\Gamma^{\,\delta-2} \,|(ad_\Gamma^{\,\delta}\,\eta_{r,n}^{\,\prime}\nabla\varphi)|
\;.
\]
Then   the Cauchy--Schwarz inequality gives
\begin{eqnarray*}
|J_{2,n}|^2&\leq& \tilde\alpha^2(\log n)^{-2}\Big(\int_{\Gamma_{\!r,n}} (ad_\Gamma^{\,\delta-2} \varphi)^2 \Big)\,
\Big(r^{-2}\int_{\Gamma_{\!r,n}} {\!\!\!}|(ad_\Gamma^{\,\delta}\,\eta_{r,n}^{\,\prime}\nabla\varphi)|^2\Big)\\[5pt]
&=&\tilde \alpha^2(\log n)^{-2}\Big(\int_{\Gamma_{\!r,n}} {\!\!\!}(ad_\Gamma^{\,\delta-2} \varphi)^2 \Big)\;J_n
\;.
\end{eqnarray*}
But if $\delta\geq 2$ then the integral is bounded by $\lambda \|\varphi\|_2^2$ or if $\delta\in[\delta_0+\delta_1r, 2\rangle$ it can be bounded by the
Rellich inequality as above.
Both cases can be handled similarly.
For example, in the latter case
\[
|J_{2,n}|^2\leq (4\tilde\alpha/\alpha_{2,r})^2(\log n)^{-2}\,\|H_{\!F}\varphi\|_2^2\;J_n
\;.
\]
Then since $J_n\leq |J_{1,n}|+|J_{2,n}|$ and $|J_{1,n}|$ is bounded by (\ref{eun3.4}) one obtains bounds 
\[
|J_{2,n}|^2\leq \sigma^2\,(\log n)^{-4}\,\|H_{\!F}\varphi\|_2^4+2\,\tau\,(\log n)^{-2}\,\|H_{\!F}\varphi\|_2^2\;|J_{2,n}|
\]
with  $\sigma,\tau>0$.
This  immediately leads to the conclusion that 
\[
|J_{2,n}|\leq 
(\sigma^2+\tau^2)^{1/2}\,(\log n)^{-2}\,\|H_{\!F}\varphi\|_2^2
\;.
\]
Alternatively, if $\delta\geq 2$ then
\[
|J_{2,n}|\leq 
(\sigma^2+\tau^2)^{1/2}\,(\log n)^{-2}\,\|\varphi\|_2\|H_{\!F}\varphi\|_2
\;.
\]
Hence $J_{2,n}\to0$ as $n\to\infty$  for all $\delta\geq \delta_0+\delta_1 r$.
Therefore  the third term on the right hand side of (\ref{eun3.1}) converges to zero for this range of $\delta$.

 It follows by combination of these observations that 
\[
\lim_{n\to\infty}\|H_{\!F}(\varphi-\varphi_{r,n})\|_2=0
\;.
\]
Therefore one concludes that $\|\varphi_{r,n}-\varphi\|_2\to0$ as $n\to\infty$ and $\|H(\varphi_{r,n}-\varphi_{r,m})\|_2\to0$ as $n,m\to\infty$.
Hence $\varphi\in D(H)$ and $D(H_{\!F})=D(H)$.
Thus $H$ is 
self-adjoint if $\delta\geq \delta_0+\delta_1 r$.
Finally remark that these arguments are valid for arbitrarily small $r$.
Therefore $H$ is self-adjoint for all $\delta>\delta_0=2-(d-d_H)/2$.
\hfill$\Box$

\bigskip

In the first case covered by Theorem~\ref{tun5.1}, the case that $\Omega$ is a $C^2$-domain, $d_H=d-1$ and $H$ is self-adjoint if $\delta>3/2$.
But if $\Omega=\Ri^d\backslash\{0\}$ then $d_H=0$ and the condition for self-adjointness is $\delta>2-d/2$.
Finally in the third case one can variously choose $\Pi$ with $d_H=1, 2, \ldots, d-1$.
Since the proof of  Theorem~\ref{tun5.1} relies completely on estimates in a boundary layer the conclusions extend to a much broader class of domains constructed using the special classes of the theorem as building blocks.
A straightforward illustration  is given by the exterior domain $\Omega=\Ri^d\backslash\Zi^d$.
Then the boundary $\Gamma$ is the lattice $\Zi^d$ and if $r$ is small the boundary layer $\Gamma_{\!\!r}$ decomposes as a countable
union of balls $B_j$ of radius $r$ centred at the points $x_j\in \Zi^d$.
But  then the theorem applied to $\Omega_j=\Ri^d\backslash\{x_j\}$ establishes  that if $\delta>2-d/2$ then $H_{\!F}\varphi_j=H\varphi_j$ for each $\varphi\in D(H_{\!F})$ with $\supp\varphi_j\subset B_j$.
Hence, by the localization results of Subsection~\ref{S2-2}, 
 $H_{\!F}\varphi=H\varphi$ for each $\varphi\in D(H_{\!F})$ with $\supp\varphi\subset \Gamma_{\!\!r}$.
 Thus $H$ is self-adjoint.
 This example is typical of models used to analyze crystalline solids.
 Of course the argument does not rely on the lattice symmetry and applies equally well to domains $\Omega=\Ri^d\backslash S$ where $S$ is a countable family of positively separated 
 points.
 Alternatively it applies if $S$ is the union of positively separated $C^2$-domains and one then deduces that the condition  $\delta>3/2$ is still sufficient for self-adjointness.

\smallskip

Nenciu and Nenciu \cite{NeN} derived conditions for self-adjointness for operators of the type we consider
by a completely different approach which makes a direct comparison  with our results rather difficult.
The major difference is that in \cite{NeN}  self-adjointness of the operators $H=-\divv(C\nabla)$ is reformulated as a problem for
Laplace--Beltrami operators on Riemannian manifolds.
This reformulation also  follows the ideas of Agmon \cite{Agm1} used in the earlier discussion of Rellich inequalities.
It consists of endowing $\Omega$ with a 
Riemannian topology corresponding to the metric defined by $C^{-1}$.
Then one can appeal to known results for diffusion operators  on manifolds.
Such results are well established if the manifold is complete but are not well understood otherwise.
The completeness or lack of completeness is directly related to boundary behaviour and 
the sufficiency criteria for self-adjointness in \cite{NeN} are given in terms of geometric
properties and the dependence on the coefficients $C$ is, in the words of Nenciu and Nenciu,  ``somewhat implicit''.
This contrasts with our approach which begins with explicit conditions on the degeneracy of the coefficients.
There are, however, several features in common between the two approaches.
The applications of the abstract results given in \cite{NeN}  appear to require the $C^2$-smoothness of the boundary,
  a condition  we have also relied upon.
Moreover, both sets of results are local and extend to domains whose boundaries are the union of positively separated parts.
It would be fruitful to further explore the possibility of combining the different arguments of the two approaches.
Next, however, we turn to a distinct problem, the derivation of necessary conditions for self-adjointness.

The derivation of  necessary conditions appears to be an underdeveloped subject.
We give an argument which starts from the earlier results \cite{LR} for $L_1$-uniqueness, i.e.\ for Markov uniqueness.
In its current form the argument only applies to $C^2$-domains and it
also relies on a slightly stronger degeneracy assumption on the derivatives of the coefficients.

\begin{thm}\label{tun5.11}
Assume  that $\Omega$ is  a $C^2$-domain.
Further assume that the coefficients $C$ of the elliptic operator $H$ satisfy
 \begin{equation}
\hspace{-13mm}\textstyle{\inf _{r\in\langle0,r_0]}}\;\textstyle{\sup_{x\in\Gamma_{\!\!r}}}\|C(x)d_\Gamma(x)^{-\delta}-a(x)I\|=0
 \label{eun5.2}
 \end{equation}
 and
\begin{equation}
\textstyle{\sup_{x\in\Gamma_{\!\!r_0}}}|(\divv (Cd^{\,-\delta})).(\nabla d_\Gamma)(x)|<\infty
\;. \label{eun5.21}
 \end{equation}

It follows that  $\delta\geq 3/2$  is necessary  for $H$ to be self-adjoint and  $\delta>3/2$ is sufficient.
\end{thm}

The  conclusion of the  condition for sufficiency, $\delta>3/2$, is restated since the assumption (\ref{eun5.21}) is slightly different to the assumption (\ref{eun3.21}) used in  Theorem~\ref{tun5.1}.
In fact the  combined assumptions (\ref{eun5.2})  and (\ref{eun5.21}) are slightly stronger than the earlier pair  of assumptions (\ref{euns3.2}) and (\ref{eun3.21}).
This follows from  the identity
\[
(\divv C).(\nabla d_\Gamma)\,d_\Gamma^{\,-\delta+1}=(\divv (Cd_\Gamma^{\,-\delta})).(\nabla d_\Gamma)\,d_\Gamma+\delta (\nabla d_\Gamma,C\nabla d_\Gamma)\,d_\Gamma^{\,-\delta}
\]
which, in combination with  (\ref{eun5.2})  and (\ref{eun5.21}), establishes that 
condition (\ref{eun3.21}) is valid.
The distinction between the current restrictions and the earlier ones is illustrated by the following example.

\begin{exam} \label{exun5.1}
Let $C=ad_\Gamma^{\,\delta}I+B_rd_\Gamma^{\,\delta+\gamma}$ where $B_r$ are positive matrices and $\gamma> 0$.
Assume $\|B_r\|$ and $|\divv B_r|$ are uniformly bounded on $\Gamma_{\!\!r}$.
Then $\|C(x)d_\Gamma(x)^{-\delta}-a(x)I\|\leq c_0d_\Gamma(x)^{\gamma}$ on $\Gamma_{\!\!r}$
and (\ref{eun5.2}) is satisfied.
Moreover, 
\[
|(\divv C).(\nabla d_\Gamma)|\,d_\Gamma^{-\delta+1}\leq c_1(1+d_\Gamma)(1+d_\Gamma^\gamma)
\;.
\]
Therefore  (\ref{eun3.21}) is also valid for all $\gamma>0$.
Nevertheless
\[
(\divv (Cd^{\,-\delta})).(\nabla d_\Gamma)=(\nabla\! a).(\nabla d_\Gamma)+(\divv B_r).(\nabla d_\Gamma)d_\Gamma^{\,\gamma}+\gamma(\nabla d_\Gamma, B_r\!\nabla d_\Gamma)d_\Gamma^{\,\gamma-1}
\]
and so (\ref{eun5.21})
 is only valid for all choices of $B_r$ if $\gamma\geq 1$.

\end{exam}

Now we turn to the proof of the theorem.
\smallskip

\noindent{\bf Proof of Theorem~\ref{tun5.11}}$\;$
The sufficiency condition follows from the foregoing discussion.
The necessary  condition $\delta\geq 3/2$ follows from  an argument by contradiction. 
It is divided into three steps.

First observe that if $\delta<2-(d-d_H)=1$ then $H$ is not Markov unique and consequently not self-adjoint.
This follows from the characterization of Markov uniqueness given by Theorem~1.1 of \cite{LR}.
This theorem is applicable since the coefficients satisfy the condition (\ref{euns2.2}).
Therefore, for self-adjointness, it is necessary that $\delta\geq 1$.
Hence  we  now assume that  $\delta\in[1,3/2\rangle$ and argue that $H$ is not self-adjoint.
The proof is divided into  two cases $ \delta\in\langle1,3/2\rangle$ and $\delta=1$ 
and for both cases  we need some extra information on Markov uniqueness.

The Friedrichs' extension $H_{\!F}$ of  $H$ is determined by the Dirichlet form $h$.
In addition there is a Markov extension $H_{\!N}$  corresponding to a Dirichlet form extension $h_N$ with the largest possible domain
\[
D(h_{N})=\{\varphi\in W^{1,2}_{\rm loc}(\Omega):\;\Gamma_{\!\!\delta}(\varphi)+\varphi^2\in L_1(\Omega)\}
\]
(see \cite{RSi4} \cite{LR} and references therein).
Then $H$ is Markov unique if and only if $h=h_{N}$.

Next  one has the following criterion for self-adjointness of $H$.

\begin{lemma}\label{l3} The following conditions are equivalent.
\begin{tabel}
\item\label{l3-1}
$H$ is self-adjoint,
\item\label{l3-2}
$D(H^*)\subseteq D(h)$,
\item\label{l3-3}
$H$ is Markov unique and $D(H^*)\subseteq D(h_{N})$.
\end{tabel}
\end{lemma}
\proof\ 
\ref{l3-1}$\Rightarrow$\ref{l3-2}$\;$
If $H$ is self-adjoint then $H^*=H=H_{\!F}$.
Hence $D(H^*)=D(H_{\!F})\subseteq D(h)$. 

\smallskip

\noindent\ref{l3-2}$\Rightarrow$\ref{l3-1}$\;$
Each self-adjoint extension $K$  of $H$ is a restriction of $H^*$.
Thus $D(K)\subseteq D(h)$. 
Hence $K=H_{\! F}$ (see \cite{Kat1}, Theorem~VI.2.11).
Therefore  $H_{\! F}$ is the unique self-adjoint extension of $H$.

 \smallskip

\noindent\ref{l3-1}$\Rightarrow$\ref{l3-3}$\;$
Self-adjointness implies Markov uniqueness, i.e.\ $h=h_{N}$.
Then Condition~\ref{l3-3} follows as before.

\smallskip

\noindent \ref{l3-3}$\Rightarrow$\ref{l3-2}$\;$
If $H$ is  Markov unique then $h=h_N$
and the implication follows immediately.
\hfill$\Box$
 
\bigskip

Now, returning to the proof of the theorem,  assume that $ \delta\in\langle1,3/2\rangle$.
Then $H$ is Markov unique by Theorem~1.1 of \cite{LR} but we next  
construct a $\nu_\delta\in D(H^*)$ such that $\nu_\delta\not\in D(h_N)$.
Thus Condition~\ref{l3-3} of  Lemma~\ref{l3} is false and consequently $H$ is not self-adjoint.
The function $\nu_\delta$ is a version of $d_\Gamma^{\,1-\delta}$ localized at a point of the boundary $\Gamma$.

Fix a $\chi\in C_c^\infty(\Ri^d)$ such that $S=\supp\chi\cap \Omega \subset\overline{\Gamma}_{\!\!r}$ has non-zero measure.
Then set $\nu_\delta=d_\Gamma^{\,1-\delta}\chi$.
Since $\delta<3/2$ it is readily verified that $\nu_\delta\in L_2(\Omega)$.
Next assume there is an $s<r$ such that $\chi=1$ on $S_s=S\cap B_s$ where $B_s$ is a ball of radius $s$ centred on $\Gamma$.
Now we argue that $\nu_\delta\not\in D(h_N)$.

It follows from the Leibniz rule and the Cauchy--Schwarz inequality that 
\begin{eqnarray*}
\Gamma(\nu_\delta)&=&  \Gamma(d_\Gamma^{\,1-\delta})\chi^2 +2\,d_\Gamma^{\,1-\delta}\, \Gamma(d_\Gamma^{\,1-\delta},\chi)\chi
+d_\Gamma^{\,2(1-\delta)}\Gamma(\chi)\\[5pt]
&\geq& (1/2) \Gamma(d_\Gamma^{\,1-\delta})\chi^2 -2\,d_\Gamma^{\,2(1-\delta)}\Gamma(\chi)
\;.
\end{eqnarray*}
Since $\Gamma(\chi)$ is bounded and  $\delta<3/2$ it follows that  the last term is integrable.
But 
\[
\Gamma(d_\Gamma^{\,1-\delta})=(\delta-1)^2d_\Gamma^{\,-2\delta}(\nabla d_\Gamma, C\nabla d_\Gamma)
\geq \sigma_{\!r}(\delta-1)^2\,ad_\Gamma^{\,-\delta}
\]
by the lower bound (\ref{euns2.2}).
Since $\delta>1$ it then follows that 
\[
 \|\Gamma(d_\Gamma^{\,1-\delta})\chi^2\|_1\geq  \sigma_{\!r}(\delta-1)^2 \int_{\!S_s}ad_\Gamma^{\,-\delta}=\infty
 \;.
  \]
 Thus $\nu_\delta\not\in D(h_N)$.
Next we prove that $\nu_\delta\in D(H^*)$.

It suffices to establish that there is a $c>0$ such that $|(H\varphi, \nu_\delta)|\leq c\,\|\varphi\|_2$
for all $\varphi\in C_c^\infty(\Gamma_{\!\!r})$.
But if $r$ is sufficiently small then $d_\Gamma$ is a $C^2$-function on $\Gamma_{\!\!r}$ (see Subsection~\ref{S2-5}).
Therefore  $\nu_\delta\in C^2_c(\Gamma_{\!\!r})$ and one calculates as before that 
\[
(H\varphi,\nu_\delta)=(\varphi, (Hd_\Gamma^{\,1-\delta})\chi)-2\,(\varphi, \Gamma(d_\Gamma^{\,1-\delta},\chi))+(\varphi, d_\Gamma^{\,1-\delta}(H\chi))
\;.
\]
Again $d_\Gamma^{\,1-\delta}(H\chi)\in L_2(\Gamma_{\!\!r})$ because $H\chi$ is bounded and $\delta<3/2$.
Moreover,
\begin{eqnarray*}
|\Gamma(d_\Gamma^{\,1-\delta},\chi)|&\leq &\Gamma(d_\Gamma^{\,1-\delta})^{1/2}\,\Gamma(\chi)^{1/2}\\[5pt]
&=&(\delta-1)\,d_\Gamma^{\,-\delta}(\nabla d_\Gamma, C\nabla d_\Gamma)^{1/2}(\nabla\chi, C\nabla\chi)^{1/2}
\leq \tau_{\!r}(\delta-1) a|\chi|^2
\end{eqnarray*}
where the last step uses (\ref{euns2.2}).
Therefore $\Gamma(d_\Gamma^{\,1-\delta},\chi)\in L_2(\Gamma_{\!\!r})$.
Next
\begin{eqnarray*}
(\delta-1)^{-1}(Hd_\Gamma^{\,1-\delta})=\divv(Cd_\Gamma^{\,-\delta}\nabla d_\Gamma)
=(\divv(Cd_\Gamma^{\,-\delta})).(\nabla d_\Gamma)+\Tr(Cd_\Gamma^{\,-\delta}D^2d_\Gamma)
\;.
\end{eqnarray*}
Hence 
\[
|(Hd_\Gamma^{\,1-\delta})\chi|\leq (\delta-1)\Big(|(\divv(Cd_\Gamma^{-\,\delta}).(\nabla d_\Gamma))|+\|Cd_\Gamma^{\,-\delta}\|\,\Tr(|D^2d_\Gamma|)\Big)|\chi|\in L_\infty(S\cap\Gamma_{\!\!r})
\]
by use of (\ref{eun2.300}), (\ref{eun5.2}) and (\ref{eun5.21}).
Combination of these estimates immediately leads to the conclusion that $\nu_\delta\in D(H^*)$.
Therefore $H$ is not self-adjoint.

Finally the case $\delta=1$ is proved similarly but one sets $\nu_1=-(\log d_\Gamma)\chi$.
\hfill$\Box$

\bigskip

In principle this method of proof could be extended to establish necessary conditions for  the other cases $\Omega=\Ri^d\backslash\{0\}$ and $\Omega=\Ri^d\backslash\overline\Pi$ covered by Theorem~\ref{tun5.1}.
 In both these cases  Markov uniqueness follows if $\delta\geq 2-(d-d_H)$ by Theorem~1.1 of \cite{LR}.
Therefore it would appear possible to repeat the foregoing arguments with $\nu_\delta=d_\Gamma^{\,2-(d-d_H)-\delta}\chi$, or $\nu_\delta=(\log d_\Gamma)\chi$
if $\delta=2-(d-d_H)$.
Then it is not difficult to verify that $\nu_\delta\not\in D(h_N)$ if $\delta<2-(d-d_H)/2$ but it is not at all clear that $\nu_\delta\in D(H^*)$.
In fact this appears to require stronger conditions on the derivatives of the coefficients than (\ref{eun3.21}) or (\ref{eun5.21}).
This is illustrated by the following example.

\begin{exam}\label{ex5.2} Let $\Omega=\Ri^d\backslash\{0\}$ and $C(x)=a(x)|x|^\delta I$ with $a\in W^{1,\infty}(\Omega)$. 
Then clearly (\ref{euns3.2}) is valid and since $a$ is bounded (\ref{eun3.21}) is equivalent to  $x\in \Gamma_{\!\!r}\mapsto |(x.\nabla a)|\in L_\infty(\Gamma_{\!\!r})$ for small $r>0$
where $\Gamma_{\!\!r}$ is now  the ball of radius $r$ centred at the origin.
It then  follows from Theorem~\ref{tun5.1} that the condition $\delta>2-d/2$ is sufficient for self-adjointness of $H$.
Now consider the proof of necessity.

If $\delta< 2-d$ then $H$ is not Markov unique by \cite{LR}. 
Hence $H$ is not self-adjoint.
Now assume $\delta> 2-d$ and 
 set $\nu_\delta(x)=|x|^{2-d-\delta}\chi(x)$ where  $\chi\in C_c^\infty({\overline\Gamma}_{\!\!r})$ with $\chi=1$ on $\Gamma_{\!\!s}$ for some $s\in\langle0,r\rangle$.
Then it follows that $\nu_\delta\in L_2(\Gamma_{\!\!r})$ if $\delta<2-d/2$
but $C|\nabla\nu_\delta|^2\not\in L_1(\Gamma_{\!\!r})$ because $\delta>2-d$.
Thus $\nu_\delta\not\in D(h_N)$ if $2-d<\delta<2-d/2$.
Now, however,
\[
(H\varphi, \nu_\delta)=(d+\delta-2)\int_{\Gamma_{\!\!r}}\varphi(x)\chi(x)(x.\nabla a)|x|^{-d}+\int_{\Gamma_{\!\!r}}\varphi(x)\,R(x)
\]
for all $\varphi\in C_c^\infty(\Omega)$ where $R$ is bounded on $\Gamma_{\!\!r}$ since it depends on derivatives of $\chi$ which are zero on $\Gamma_{\!\!s}$.
Thus $\nu_\delta\in D(H^*)$ if and only if  $x\mapsto|(x.\nabla a)||x|^{-d}\in L_2(\Gamma_{\!\!r})$.
(A similar conclusion follows in the case $\delta=2-d$ with $\nu_\delta(x)=(\log|x|)\chi(x)$.)
But in this example (\ref{eun3.21}) is equivalent to boundedness of $x\in\Gamma_{\!\!r}\mapsto|(x.\nabla a)|\in L_\infty(\Gamma_{\!\!r})$
and (\ref{eun5.21}) is equivalent to $x\in\Gamma_{\!\!r}\mapsto|(x.\nabla a)||x|^{-1}\in L_\infty(\Gamma_{\!\!r})$.
Both the $L_\infty$-conditions follow from the assumption $a\in W^{1,\infty}(\Omega)$ but the $L_2$-condition is generally much stronger.
Nevertheless one does conclude that if $x\mapsto|(x.\nabla a)||x|^{-d}\in L_2(\Gamma_{\!\!r})$
then the condition  $\delta\geq 2-d/2$ is necessary for self-adjointness of $H$.

Finally note that if $a$ is a radial function then all these conditions simplify since $|(x.\nabla a)||x|^{-1}=a^\prime(|x|)$.
But in this case one can completely analyze the example by passing to radial coordinates and using the Weyl limit point theory.
One then has self-adjointness of $H$ if and only if $\delta\geq2-d/2$.
If, however,  $a$ is an angular function, i.e.\ a function on the unit sphere, then this approach is not possible.
Nevertheless one then has  $|(x.\nabla a)|=0$ and the foregoing arguments all apply.
\end{exam}

It is of course possible that the foregoing arguments apply to $\Ri^d\backslash\{0\}$ or $\Ri^d\backslash\overline\Pi$ with a different choice of $\nu_\delta$ although the singular
factor $d_\Gamma^{\,2-(d-d_H)-\delta}$ seems compelling.

\section{Concluding remarks}\label{S4}

In conclusion we discuss some potential improvements in the foregoing results, some apparent difficulties in their extension and some particular examples.

First note that Theorem~\ref{tun5.1}  
 leaves open the question whether the critical relation $\delta=2-(d-d_H)/2$ is sufficient  for self-adjointness
 although Theorem~\ref{tun5.11} establishes that this condition is necessary in the $C^2$-case.
The current  sufficiency arguments fail  since the Rellich inequality gives no information at the critical value.
For example, the inequality  for the Laplacian on $\Ri^d\backslash\{0\}$ is non-trivial if and only if $d\geq 5$.
One possible way of circumventing this difficulty might be with  a modified Rellich inequality containing a logarithmic factor, an inequality of the type given
by Proposition~\ref{pex2.1} but with an additional factor $(\log d_\Gamma)^{-2}$ in the right hand side.
Critical Hardy inequalities for $\Omega=\Ri^d\backslash\{0\}$ with a logarithmic modification were considered by  Solomyak \cite{Sol} and recently have been analyzed intensively (see, for example, \cite{HoK} \cite{IoM}  \cite{Taka} \cite{SaT} and references therein).
It is unclear whether one might expect similar weighted Rellich inequalities for $\Ri^d\backslash\{0\}$ or more general~$\Omega$.

Secondly, it is possible that modification of the degeneracy conditions (\ref{euns3.2}) and (\ref{eun3.21}) might lead to improved results.
The one-dimensional case is well-understood
(see \cite{CMP}, Proposition~3.5, for  finite intervals $\Omega$ or \cite{RSi3}, Theorem~2.4, for semi-infinite intervals).
Then the operator $H$ has the action $H\varphi=-(c\,\varphi')'$.
Now suppose for simplicity that $\Omega=\langle0,\infty\rangle$.
Then it follows from Theorem~2.4 of \cite{RSi3} that 
$H$ is essentially self-adjoint on $C_c^\infty(0,\infty)$ if and only if  $x>0\mapsto\nu(x)=\int^r_xds\,c(s)^{-1}\not\in L_2(0,r)$ for some small $r\in\langle0,1]$.
The result leads to a  simplified one-dimensional version of  Theorems~\ref{tun5.1}  and \ref{tun5.11}.
\begin{obs} \label{O4.1}Assume that $\limsup_{x\to0}c(x)x^{-\delta}=a>0$ then $H$ is essentially self-adjoint  on $C_c^\infty(0,\infty)$ if and only if $\delta\geq3/2$.
\end{obs}
\proof\ Choose $r>0$ such that $\sup_{\{x\in\langle0,r]\}}|c(x)x^{-\delta}-a|<a/2$.
Then it follows that $|c(x)^{-1}-(ax^\delta)^{-1}|<(2c(x))^{-1}$. 
Hence $(2/3)(ax^\delta)^{-1}\leq c(x)^{-1}\leq 2(ax^\delta)^{-1}$.
Therefore the condition $\nu\not\in L_2(0,r)$ is equivalent to $\int^r_0dx \,x^{2(1-\delta)}=\infty$ or,  equivalently, $\delta\geq 3/2$.
\hfill$\Box$

\bigskip

A similar conclusion is valid on a finite interval and the values of $a$ and $\delta$ may differ from endpoint to endpoint.
This result is notable since it does not require any explicit bound on the derivative $c'$ of the coefficient $c$. 
It indicates that there could be variants of the multi-dimensional results which  do not require the boundedness condition 
(\ref{eun3.21}) on $\divv C$. 

Thirdly, it is not clear whether $C^2$-regularity of the boundary is essential for the self-adjointness results of Section~\ref{S3}.
It is possible the conclusions are valid for $C^{1,1}$-boundaries.
It is known that  $\Omega$ is a $C^{1,1}$-domain if and only if it  satisfies  both an interior and exterior sphere condition (An extensive discussion can be found in  \cite{Barb}).
Thus there are internal and external boundary layers in which the  distance $d_\Gamma$ is differentiable. 
Moreover the signed distance  is a $C^{1,1}$-function.
(See \cite{DeZ} Section~7.8 for a complete analysis and background references.)
The essential estimate on the distance function  in the discussion of $C^2$-domains was the bound $|(d_\Gamma\nabla^2d_\Gamma)(x)|\leq \gamma r$ for all $x\in\Gamma_{\!\!r}$.
This was a direct consequence of the Gilbarg--Trudinger result \cite{GT} Lemma~14.17.
Although the proof of the  latter does not extend to the $C^{1,1}$-case it is likely that the basic estimate is still valid.

Fourthly, an alternative possible approach might be  by capacity estimates.
The earlier analysis of $L_1$-uniqueness in \cite{RSi4} and \cite{LR} was based on capacity arguments
with the capacity $\capp_h$ defined in terms of the quadratic form $h$.
The $L_1$-analysis began with Theorems~1.2 and 1.3 of \cite{RSi4} which
 established that $L_1$-uniqueness is equivalent both to Markov uniqueness and 
 to the capacity  condition $\capp_h(\Gamma)=0$.
A similar analysis of $L_2$-uniqueness could  then be  based on the capacity $\capp_H$ defined in terms of the operator $H$
with the aim to prove that the uniqueness property is now equivalent to $\capp_H(\Gamma)=0$.
Explicitly   the capacity  of the measurable subset  $A\subset \overline\Omega$ is given  by
 \begin{eqnarray*}
\capp_H(A)=\inf\Big\{\;\|\psi\|_{D(H_{\!F})} &&\;: \;\psi\in D(H_{\!F}) \mbox{ and  there exists   an open set  }\nonumber\\[-5pt]
&& U\subset \Ri^d
\mbox{ such that } U\supseteq A  \mbox{ and }
\psi=1  \mbox{ a.\ e. on } U\cap\Omega\;\Big\}
\;.
\end{eqnarray*}
One can apply this algorithm to calculate the capacity of bounded subsets $A$ of the boundary $\Gamma$.
Since $D(H_{\!F})\supset W^{2,2}(\Omega) $ the capacity of each $A$ is finite.
But using the fact that  $\delta\geq 2-(d-d_H)/2$ one can in fact establish that $\capp_H(A)=0$ for each bounded $A\subset \Gamma$.
The calculation follows similar lines to the proof of Proposition~3.1 in \cite{LR} but the $\eta_{r,n}$ in that proof are replaced by
$\eta_n\circ (r^{-1}d_A)$ where the $\eta_n$ are now the approximate identity constructed in Subsection~\ref{S2-3}.
This substitution results in the critical value $\delta_c=2-(d-d_H)$ in Proposition~3.1 of  \cite{LR}  being replaced by  $2-(d-d_H)/2$.
Then by general properties of the capacity one should have $\capp_H(\Gamma)=0$.
Thus if one hopes to follow the reasoning for $L_1$-uniqueness it remains to prove that $\capp_H(\Gamma)=0$ implies that $H$ is self-adjoint, i.e.\ $L_2$-unique.
A result of this nature can be established for the Laplacian on certain exterior domains
(see, for example, \cite{AdH} Corollary~5.1.15) but the  difficulties  of dealing with variable coefficients appear considerable.

\end{document}